\documentclass[11pt]{amsart}

\usepackage{hyperref}
\hypersetup{nesting=true,debug=true,naturalnames=true}
\usepackage{graphicx,amssymb,upref}
\usepackage[usenames,dvipsnames]{pstricks}
\usepackage{epsfig}
\usepackage{color}
\usepackage[utf 8]{inputenc}
\usepackage[T1]{fontenc}
\usepackage{amsmath, amsthm}
\usepackage[english]{babel}
\usepackage{amssymb}
\usepackage{latexsym}
\usepackage{graphicx,amssymb,upref,mathtools}
\usepackage[all]{xy}
\usepackage{tikz}
\usepackage{yhmath}
\usetikzlibrary{matrix} 
\usepackage{tikz-cd}

\allowdisplaybreaks

\date{\today}

\newtheorem{theorem}{Theorem}[section]

\newtheorem{corollary}[theorem]{Corollary}
\theoremstyle{definition}
\newtheorem{definition}[theorem]{Definition}
\newtheorem{example}[theorem]{Example}

\newtheorem{remark}[theorem]{Remark}

\DeclareMathOperator{\Hom}{Hom}

\newcommand{\wh}{\widehat{H}}

\title[About Hopf braces and crossed products]{About Hopf braces and crossed products} 

\begin{document}
	
\maketitle
	
\begin{center}
{\bf Ram\'on
Gonz\'{a}lez Rodr\'{\i}guez$^{1}$ and Brais Ramos P\'erez$^{2}$}.
\end{center}
	
\vspace{0.4cm}
	
\begin{center}
{\small $^{1}$ [https://orcid.org/0000-0003-3061-6685].}
\end{center}
\begin{center}	{\small  CITMAga, 15782 Santiago de Compostela, Spain.}
\end{center}
\begin{center}
{\small  Universidade de Vigo, Departamento de Matem\'{a}tica Aplicada II, E. E. Telecomunicaci\'on,
E-36310 Vigo, Spain.
\\email: rgon@dma.uvigo.es}
\end{center}
\vspace{0.2cm}

\begin{center}
{\small $^{2}$ [https://orcid.org/0009-0006-3912-4483].}
\end{center}
\begin{center}
{\small  CITMAga, 15782 Santiago de Compostela, Spain. \\}
\end{center}
\begin{center}
{\small  Universidade de Santiago de Compostela. Departamento de Matem\'aticas,  Facultade de Matem\'aticas, E-15771 Santiago de Compostela, Spain. 
\\email: braisramos.perez@usc.es}
\end{center}
\vspace{0.2cm}

\begin{abstract}  
The present article represents a step forward in the study of the following problem: If $\mathbb{A}=(A_{1},A_{2})$ and $\mathbb{H}=(H_{1},H_{2})$ are Hopf braces in a symmetric monoidal category {\sf C} such that $(A_{1},H_{1})$ and $(A_{2},H_{2})$ are matched pairs of Hopf algebras, then we want to know under what conditions the pair $(A_{1}\bowtie H_{1},A_{2}\bowtie H_{2})$ constitutes a new Hopf brace. We find such conditions for the pairs $(A_{1}\otimes H_{1},A_{2}\bowtie H_{2})$ and $(A_{1}\bowtie H_{1},A_{2}\sharp H_{2})$ to be Hopf braces, which are particular situations of the general problem described above, and we apply these results to study when the Drinfeld's Double gives rise to a Hopf brace.
\end{abstract} 

\vspace{0.2cm}

{\sc Keywords}: Symmetric monoidal category, Hopf algebra, matched pair, bicrossed product, smash product, Hopf brace, Drinfeld's Double.

{\sc MSC2020}: 18M05, 16T05, 16S35, 16S40.

\section*{Introduction}

A Hopf brace is a recent mathematical object which was born as the quantum version of a skew brace. As a generalization of W. Rump's braces (see \cite{Rump}), skew braces appeared in \cite{GV} in order to study non-degenerate solutions of the Quantum Yang-Baxter Equation and they consist of a pair of groups, $(G,\cdot)$ and $(G,\ast)$, satisfying the following compatibility condition:
\begin{gather*}
g\ast(h\cdot t)= (g\ast h)\cdot g^{-1}\cdot (g\ast t)
\end{gather*}
for all $g,h,t\in G$, where $g^{-1}$ denotes the inverse of $g$ with respect to $(G,\cdot)$. Then, a direct Hopf-theoretical generalization of skew braces is the notion of Hopf brace introduced by I. Angiono, C. Galindo and L. Vendramin in \cite{AGV} as a pair of Hopf algebras, $H_{1}=(H,1_{\circ},\circ,\epsilon,\Delta,S_{\circ})$ and $H_{2}=(H,1_{\star},\star,\epsilon,\Delta,S_{\star})$, such that the following equality holds:
\begin{gather*}\label{HBr_compat_elements}
g\star(h\circ t)=(g_{1}\star h)\circ S_{\circ}(g_{2})\circ(g_{3}\star t)
\end{gather*}
for all $g,h,t\in H$. In this introduction, Hopf braces will be denoted by $\mathbb{H}=(H_{1},H_{2})$, or by simply a pair $(H_{1},H_{2})$. It is important to note that both Hopf algebras involved in a Hopf brace share the same underlying coalgebra structure, and therefore, a Hopf brace is said to be cocommutative when both, $H_{1}$ and $H_{2}$, are cocommutative Hopf algebras. This subclass of Hopf braces is of significant importance because it is proved in \cite[Corollary 2.4]{AGV} that they also induce solutions to the Quantum Yang-Baxter equation. 

A well-studied problem in the theory of Hopf algebras (and, of course, in the theory of groups) is the factorization: Given a Hopf algebra $X$ in the category of vector spaces over a field ${\mathbb F}$, what are the conditions under which there exist $A$ and $H$, sub-Hopf algebras of $X$, such that $X=AH$? In the solution of this problem, the notion of matched pair of Hopf algebras plays a crucial role. Matched pairs were introduced in the group setting by M. Takeuchi in \cite{Tak81} and extended to Hopf algebras by S. Majid in \cite{MAJ_bicros}. So, in the Hopf algebra setting, a matched pair is a 4-tuple $(A,H,\varphi_{A},\phi_{H})$, where $A$ and $H$ are Hopf algebras, $A$ is a left $H$-module coalgebra with the left action $\varphi_{A}\colon H\otimes A\rightarrow A$ and $H$ is a right $A$-module coalgebra with the right action $\phi_{H}\colon H\otimes A\rightarrow H$, satisfying the following conditions:
\begin{gather*}
\varphi_{A}(h\otimes 1_{A})=\epsilon_{H}(h)1_{A}, \quad \varphi_{A}(h\otimes ab)=\varphi_{A}(h_{1}\otimes a_{1})\varphi_{A}(\phi_{H}(h_{2}\otimes a_{2})\otimes b),\\
\phi_{H}(1_{H}\otimes a)=\epsilon_{A}(a)1_{H},\quad \phi_{H}(hg\otimes a)=\phi_{H}(h\otimes \varphi_{A}(g_{1}\otimes a_{1}))\phi_{H}(g_{2}\otimes a_{2}),\\
\phi_{H}(h_{2}\otimes a_{2})\otimes \varphi_{A}(h_{1}\otimes a_{1})=\phi_{H}(h_{1}\otimes a_{1})\otimes \varphi_{A}(h_{2}\otimes a_{2}),
\end{gather*}
for all $a,b\in A$ and $h,g\in H$.  In \cite[Proposition 3.12]{MAJ_bicros}, S. Majid asserts that every matched pair of Hopf algebras $(A,H,\varphi_{A},\phi_{H})$ determines a new Hopf algebra structure over $A\otimes H$, the so-called bicrossed product Hopf algebra $A\bowtie H$, which is a Hopf algebra with the tensor product unit, counit and coproduct, and whose product and antipode are defined as 
\begin{gather*}
(a\otimes h)(b\otimes g)=a\varphi_{A}(h_{1}\otimes b_{1})\otimes \phi_{H}(h_{2}\otimes b_{2})g,\\
S_{A\bowtie H}(a\otimes h)=\varphi_{A}(S_{H}(h)_{1}\otimes S_{A}(a)_{1})\otimes \phi_{H}(S_{H}(h)_{2}\otimes S_{A}(a)_{2})
\end{gather*}
for all $a,b\in A$ and $h,g\in H$. Hence, matched pairs solve the problem of factorization in the Hopf algebra setting as was proved in \cite[Theorem 7.2.3]{MAJ_Book}: $X$ admits an exact factorization of the form $X=AH$, where $A$ and $H$ are sub-Hopf algebras of $X$, if and only if there exists a matched pair $(A,H,\varphi_{A}, \phi_{H})$ such that $A\bowtie H$ is a Hopf algebra isomorphic to $X$. This problem and the construction of the bicrossed product Hopf algebra from a matched pair have been studied in other contexts such as weak Hopf algebras (see \cite{BohmTorre}), Hopf quasigroups (see \cite{RGON2}) or quasigrupoids (see \cite{RGON3}).

Based on the previous background, why not studying matched pairs and the factorization problem for Hopf braces? In order to tackle this problem, the first step is to define the notion of matched pair of Hopf braces or, in other words, given Hopf braces $\mathbb{A}=(A_{1},A_{2})$ and $\mathbb{H}=(H_{1},H_{2})$ such that $(A_{1},H_{1},\varphi_{A}^1,\phi_{H}^1)$ and $(A_{2},H_{2},\varphi_{A}^2,\phi_{H}^2)$ are matched pairs of Hopf algebras, to find out the conditions under which the pair of Hopf algebras $(A_{1}\bowtie H_{1},A_{2}\bowtie H_{2})$ is a Hopf brace. In general, obtaining such conditions is not a simple problem. The first contributions were made by A. Agore in \cite{AGORE}, who combines a Hopf algebra and a Hopf brace by using crossed products in order to construct new Hopf braces where one of the Hopf algebras involved is always the bicrossed product Hopf algebra described earlier. In this paper, despite not having achieved a complete solution to this problem, the desired conditions have been found in two particular cases:
\begin{itemize}
\item[(i)] In Theorem \ref{mainth2}, the conditions for $(A_{1}\otimes H_{1},A_{2}\bowtie H_{2})$ to be a Hopf brace have been found, which is a particular situation of the general problem by considering that the actions $\varphi_{A}^1$ and $\phi_{H}^1$ are the trivial actions.
\item[(ii)] Moreover, in Theorem \ref{mainth}, the conditions for $(A_{1}\bowtie H_{1},A_{2}\sharp H_{2})$ to be a Hopf brace have been deduced, where $\sharp$ denotes the well-known smash product Hopf algebra. Again, this is a particular case of the general problem by assuming that $\phi_{H}^2$ is the trivial action.
\end{itemize}

In \cite{TD}, given a finite dimensional Hopf algebra $H$, Y. Doi and M. Takeuchi proved that the Drinfeld's Double Hopf algebra $D(H)$ can be obtained as the bicrossed product Hopf algebra determined by a matched pair $(A,H,\varphi_{A}, \phi_{H})$, where $A=(H^{cop})^{\ast}$ and the actions $\varphi_{A}$ and $\phi_{H}$ are defined in a particular way that will be made explicit in the development of the article. The Drinfeld's Double construction is a classical structure in the theory of Hopf algebras introduced by V. Drinfeld in \cite{Drin} whose importance lies in the fact that $D(H)$ is a quasitriangular Hopf algebra (independently of the fact that the initial Hopf algebra $H$ is not), and thus, it induces solutions to the Quantum Yang-Baxter Equation. Since $D(H)=(H^{cop})^{\ast}\bowtie H$ as a consequence of the result by Y. Doi and M. Takeuchi, the main theorems mentioned above can be applied to this situation in order to obtain conditions under which the Drinfeld's Double is one of the Hopf algebras that is part of a Hopf brace structure, a situation that has not been studied so far. To be precise, if $T(H)\coloneqq (H^{cop})^{\ast}\otimes H$ the tensor product Hopf algebra, the following results are obtained:
\begin{itemize}
\item[(i)] On the one hand, in Theorem \ref{D(H)_cocommutative_HBr} it is proved that $(T(H),D(H))$ is always a Hopf brace in case $H$ is a cocommutative Hopf algebra. Moreover, in the non-cocommutative setting, Corollary \ref{D(H)_cocomclass_HBr} sets the conditions under which $(T(H),D(H))$ is a Hopf brace when $((H^{cop})^{\ast},\varphi_{(H^{cop})^{\ast}})$ is in the cocommutativity class of $H$, which is a weaker assumption than cocommutativity.
\item[(ii)] On the other hand, Theorem \ref{D(H)_commutative_HBr} asserts that $(D(H),T(H))$ is always a Hopf brace under commutativity of $H$. However, without making any assumptions, $(D(H),T(H))$ is a Hopf brace if and only if equality \eqref{D(H)_HBr_nocommutativity} holds, as can be consulted in Corollary \ref{D(H)_noncomm_HBr}.
\end{itemize} 

\section{Preliminaries}
Throughout  this paper, we will be working on a strict symmetric monoidal setting. Thus, for us {\sf C} will denote a strict symmetric monoidal category with tensor product $\otimes$, unit object $K$ and natural isomorphism of symmetry $c$.
Along the paper we will use the notation $P\otimes f$ and $f\otimes P$ to refer to $id_{P}\otimes f$ and $f\otimes id_{P}$, respectively, for all $f\colon M\rightarrow N$ morphism in {\sf C} and for all $P$ object in ${\sf C}$. 

The following concepts will be used repeatedly in the development of the article.
\begin{definition}
A triple $A=(A,\eta_{A},\mu_{A})$ is said to be an algebra in {\sf C} if $\eta_{A}\colon K\rightarrow A$ (unit) and $\mu_{A}\colon A\otimes A\rightarrow A$ (product) are morphisms in {\sf C} such that  following conditions hold:
\begin{gather*}
\mu_{A}\circ (\eta_{A}\otimes A)=id_{A}=\mu_{A}\circ(A\otimes\eta_{A}),\quad\mu_{A}\circ(\mu_{A}\otimes A)=\mu_{A}\circ(A\otimes\mu_{A}).
\end{gather*}

 Note that, given $B=(B,\eta_{B},\mu_{B})$ another algebra in {\sf C}, the tensor product $A\otimes B$ admits an algebra structure whose unit and product are defined as follows:
$$
\eta_{A\otimes B}\coloneqq \eta_{A}\otimes\eta_{B}, \quad
\mu_{A\otimes B}\coloneqq (\mu_{A}\otimes\mu_{B})\circ(A\otimes c_{B,A}\otimes B).
$$

Dually, a triple $C=(C,\varepsilon_{C},\delta_{C})$ is a coalgebra in {\sf C} if  $\varepsilon_{C}\colon C\rightarrow K$ (counit) and $\delta_{C}\colon C\rightarrow C\otimes C$ (coproduct) are morphisms in {\sf C} satisfying the following conditions:
\begin{gather*}
(\varepsilon_{C}\otimes C)\circ\delta_{C}=id_{C}=(C\otimes\varepsilon_{C})\circ\delta_{C},\quad(\delta_{C}\otimes C)\circ\delta_{C}=(C\otimes\delta_{C})\circ\delta_{C}.
\end{gather*}

 Likewise, if $D=(D,\varepsilon_{D},\delta_{D})$ is another coalgebra in {\sf C}, the tensor product $C\otimes D$ is a coalgebra with the following counit and coproduct:
$$
\varepsilon_{C\otimes D}\coloneqq\varepsilon_{C}\otimes\varepsilon_{D},\quad
\delta_{C\otimes D}\coloneqq(C\otimes c_{C,D}\otimes D)\circ(\delta_{C}\otimes\delta_{D}).
$$
\end{definition}

\begin{definition}
Let $A$ and $B$ be algebras in {\sf C}. A morphism $f\colon A\rightarrow B$ in {\sf C} is said to be an algebra morphism if the following equalities hold:
$$
f\circ\eta_{A}=\eta_{B}, \quad f\circ\mu_{A}=\mu_{B}\circ(f\otimes f).
$$

\noindent In other words, $f$ is an algebra morphism if $f$ preserves the unit and it is multiplicative.

Similarly, given $C$ and $D$ coalgebras in {\sf C}, a morphism $g\colon C\rightarrow D$ is a coalgebra morphism if the conditions 
$$
\varepsilon_{D}\circ g=\varepsilon_{C}, \quad\delta_{D}\circ g=(g\otimes g)\circ\delta_{C}
$$
hold, i.e., $g$ preserves  the counit and it is comultiplicative.
\end{definition}

Let $C$ be a coalgebra and let $A$ be an algebra in {\sf C}. The set of morphisms in {\sf C} between $C$ and $A$, denoted by $\operatorname{Hom}_{{\sf C}}(C,A)$, is a monoid with the convolution product, which is defined as follows:
\[f\ast g\coloneqq\mu_{A}\circ (f\otimes g)\circ\delta_{C}.\]

The unit for this operation is given by $\eta_{A}\circ\varepsilon_{C}=\varepsilon_{C}\otimes\eta_{A}$.

\begin{definition}
A 5-tuple $B=(B,\eta_{B},\mu_{B},\varepsilon_{B},\delta_{B})$ is said to be a bialgebra in {\sf C} if $(B,\eta_{B},\mu_{B})$ is an algebra in {\sf C}, $(B,\varepsilon_{B},\delta_{B})$ is a coalgebra in {\sf C} and $\eta_{B}$ and $\mu_{B}$ are coalgebra morphisms.
\end{definition}
\begin{remark}
Note that the fact that $\eta_{B}$ and $\mu_{B}$ are coalgebra morphisms is equivalent to the fact that $\varepsilon_{B}$ and $\delta_{B}$ are algebra morphisms.
\end{remark}

\begin{definition}
Let $B$ and $D$ be bialgebras in {\sf C}. A morphism $f\colon B\rightarrow D$ in {\sf C} is a bialgebra morphism when $f$ is an algebra-coalgebra morphism.
\end{definition}

In the following definition we introduce the notion of Hopf algebra.

\begin{definition}
A 6-tuple $H=(H,\eta_{H},\mu_{H},\varepsilon_{H},\delta_{H},\lambda_{H})$ is said to be a Hopf algebra in {\sf C} when $(H,\eta_{H},\mu_{H},\varepsilon_{H},\delta_{H})$ is a bialgebra in {\sf C} such that there exists an endomorphism $\lambda_{H}\colon H\rightarrow H$, called the antipode, satisfying the following identity:
\begin{gather}\label{antipode}
\lambda_{H}\ast id_{H}=\varepsilon_{H}\otimes\eta_{H}=id_{H}\ast \lambda_{H}.
\end{gather}\end{definition}
The equation \eqref{antipode} means that the antipode, $\lambda_{H}$, is the inverse of the identity for the convolution product in $\operatorname{Hom}_{{\sf C}}(H,H)$ and, as a consequence, it is unique. In addition, other relevant properties that $\lambda_{H}$ satisfies are the following: It is antimultiplicative and anticomultiplicative, which means that the identities
\begin{gather}
\label{a-antip1} \lambda_{H}\circ\mu_{H}=\mu_{H}\circ(\lambda_{H}\otimes\lambda_{H})\circ c_{H,H},\\\label{a-antip2} \delta_{H}\circ\lambda_{H}=c_{H,H}\circ(\lambda_{H}\otimes\lambda_{H})\circ\delta_{H}
\end{gather}
hold, and also $\lambda_{H}$ preserves the unit and the counit:
\begin{gather}
\label{u-antip1} \lambda_{H}\circ\eta_{H}=\eta_{H},\\\label{u-antip2}\varepsilon_{H}\circ\lambda_{H}=\varepsilon_{H}.
\end{gather}

We will say that a Hopf algebra $H$ is commutative if $\mu_{H}\circ c_{H,H}=\mu_{H}$, and cocommutative when $c_{H,H}\circ\delta_{H}=\delta_{H}$. It is a direct consequence of identities \eqref{a-antip1} and \eqref{u-antip1} that $\lambda_{H}$ is an algebra morphism when $H$ is commutative and, under cocommutativity hypothesis, $\lambda_{H}$ is a coalgebra morphism thanks to \eqref{a-antip2} and \eqref{u-antip2}. In both situations it also holds that $\lambda_{H}^{2}=id_{H},$ and then $\lambda_{H}^{-1}=\lambda_{H}$.

Note that if $H$ is a Hopf algebra in {\sf C} whose antipode $\lambda_{H}$ is an isomorphism, then 
\begin{gather*}
H^{op}=(H,\eta_{H},\mu_{H}\circ c_{H,H},\varepsilon_{H},\delta_{H},\lambda_{H}^{-1})\text{ and }H^{cop}=(H,\eta_{H},\mu_{H},\varepsilon_{H},c_{H,H}\circ\delta_{H},\lambda_{H}^{-1})
\end{gather*}
are also Hopf algebras in {\sf C}. It results clear that, if $H$ is commutative, then $H^{op}=H$, whereas if $H$ is cocommutative, then $H^{cop}=H$.

\begin{definition}
Let $H$ and $L$ be Hopf algebras in {\sf C}. A morphism $f\colon H\rightarrow L$ in {\sf C} is a morphism of Hopf algebras if it is a bialgebra morphism, that is to say, if $f$ is simultaneously an algebra and a coalgebra morphism.
\end{definition}
It is straightforward to show that $\lambda_{L}\circ f$ and $f\circ \lambda_{H}$ are both inverses of $f$ for the convolution product in $\operatorname{Hom}_{{\sf C}}(H,L)$. As a consequence of the uniqueness of the inverse for the convolution, the following equality is obtained:
\begin{gather*}\label{morant}
f\circ\lambda_{H}=\lambda_{L}\circ f.
\end{gather*}

Once the concept of Hopf algebra has been introduced, now we define what a module over a Hopf algebra is.

\begin{definition}
Let $H$ be a Hopf algebra in {\sf C}. Given an object $M$ in {\sf C}, a pair $(M,\varphi_{M})$ is said to be a left module over $H$, where $\varphi_{M}\colon H\otimes M\rightarrow M$ is a morphism in {\sf C} called the action, if the following conditions hold:
$$
\varphi_{M}\circ(\eta_{H}\otimes M)=id_{M},\quad\varphi_{M}\circ(H\otimes\varphi_{M})=\varphi_{M}\circ(\mu_{H}\otimes M).
$$

Let $(N,\varphi_{N})$ be another left module over the Hopf algebra $H$. A morphism $f\colon M\rightarrow N$ in {\sf C} is a morphism of left $H$-modules (in other words, $f$ is $H$-linear) if the equality
\begin{gather*}
f\circ\varphi_{M}=\varphi_{N}\circ(H\otimes f)
\end{gather*}
holds.
\end{definition}

Suppose now that $(M,\varphi_{M})$ is a left $H$-module, but $M$ also has another structures, for example, we can suppose that $M$ is an algebra or a coalgebra in {\sf C}. This situation gives rise to the notions of left module algebra and left module coalgebra, which set certain compatibility conditions between the (co)product and the (co)unit with the action of the module.

\begin{definition}
Let $H$ be a Hopf algebra and $A$ an algebra in {\sf C}. It is said that $(A,\varphi_{A})$ is a left $H$-module algebra if $(A,\varphi_{A})$ is a left $H$-module such that $\eta_{A}$ and $\mu_{A}$ are morphisms of left $H$-modules, which translates into the fact that the following identities are verified:
\begin{gather*}
\label{Hmodalg1}\varphi_{A}\circ(H\otimes\eta_{A})=\varepsilon_{H}\otimes\eta_{A},\\\label{Hmodalg2}\varphi_{A}\circ(H\otimes\mu_{A})=\mu_{A}\circ\varphi_{A\otimes A},
\end{gather*} 
where $\varphi_{A\otimes A}\coloneqq (\varphi_{A}\otimes\varphi_{A})\circ(H\otimes c_{H,A}\otimes A)\circ(\delta_{H}\otimes A\otimes A)$ is the left $H$-module structure over $A\otimes A$.

On the other hand, let $C$ be a coalgebra in {\sf C}. It is said that $(C,\varphi_{C})$ is a left $H$-module coalgebra if $(C,\varphi_{C})$ is a left $H$-module such that $\varepsilon_{C}$ and $\delta_{C}$ are morphisms of left $H$-modules, which means that the equalities
\begin{gather}
\label{Hmodcoalg1}\varepsilon_{C}\circ\varphi_{C}=\varepsilon_{H}\otimes\varepsilon_{C},\\\label{Hmodcoalg2}\delta_{C}\circ\varphi_{C}=\varphi_{C\otimes C}\circ(H\otimes\delta_{C})
\end{gather} 
hold. Note that \eqref{Hmodcoalg1} and \eqref{Hmodcoalg2} are equivalent to the fact that $\varphi_{C}$ is a coalgebra morphism.
\end{definition}

\begin{remark}
Following the same line of the previous definitions it is possible to define the notions of right module and right module (co)algebra over a Hopf algebra.
\end{remark}

\begin{example}
Let $H$ be a Hopf algebra in {\sf C}. $(H,\varphi_{H}^{ad})$ is a left $H$-module algebra, where $$\varphi_{H}^{ad}\coloneqq \mu_{H}\circ(\mu_{H}\otimes\lambda_{H})\circ(H\otimes c_{H,H})\circ(\delta_{H}\otimes H),$$
the so-called adjoint action of $H$. Moreover, if $H$ is cocommutative, then $(H,\varphi_{H}^{ad})$ is a left $H$-module coalgebra too.
\end{example}

In this paper the construction of Hopf algebras using the notion of matched pair will be relevant. In our setting, the definition of matched pair of Hopf algebras is the following:
 
\begin{definition}\label{MPdef}
Let $H$ and $A$ be Hopf algebras in {\sf C}. A 4-tuple $(A,H,\varphi_{A},\phi_{H})$ is said to be a matched pair of Hopf algebras if the following conditions hold:
\begin{itemize}
\item[(i)] $(A,\varphi_{A})$ is a left $H$-module coalgebra and $(H,\phi_{H})$ is a right $A$-module coalgebra,
\item[(ii)] $\varphi_{A}\circ (H\otimes\eta_{A})=\varepsilon_{H}\otimes\eta_{A}$, i.e., $\eta_{A}$ is a morphism of left $H$-modules, 
\item[(iii)] $\phi_{H}\circ(\eta_{H}\otimes A)=\eta_{H}\otimes\varepsilon_{A}$, i.e., $\eta_{H}$ is a morphism of right $A$-modules,
\item[(iv)] $\varphi_{A}\circ(H\otimes\mu_{A})=\mu_{A}\circ (A\otimes\varphi_{A})\circ (\Psi\otimes A)$,
\item[(v)] $\phi_{H}\circ(\mu_{H}\otimes A)=\mu_{H}\circ(\phi_{H}\otimes H)\circ(H\otimes \Psi)$,
\item[(vi)] $c_{A,H}\circ\Psi=(\phi_{H}\otimes\varphi_{A})\circ (H\otimes c_{H,A}\otimes A)\circ(\delta_{H}\otimes\delta_{A})$,
\end{itemize}
where \begin{equation}\label{interw}\Psi\coloneqq (\varphi_{A}\otimes\phi_{H})\circ (H\otimes c_{H,A}\otimes A)\circ(\delta_{H}\otimes\delta_{A}).\end{equation}
\end{definition}

It is well known that every matched pair of Hopf algebras gives rise to a new Hopf algebra structure over the tensor product $A\otimes H$ in the following way: If $(A,H,\varphi_{A},\phi_{H})$ is  a matched pair of Hopf algebras, then we have that 
\[A \bowtie H=(A\otimes H,\eta_{A\triangleright\hspace{-0.04 cm}\triangleleft H},\mu_{A\triangleright\hspace{-0.04 cm}\triangleleft H},\varepsilon_{A\triangleright\hspace{-0.04 cm}\triangleleft H},\delta_{A\triangleright\hspace{-0.04 cm}\triangleleft H},\lambda_{A\triangleright\hspace{-0.04 cm}\triangleleft H}),\]
where
\begin{gather*}
\eta_{A\bowtie H}\coloneqq\eta_{A\otimes H}, \;\;\; \mu_{A\bowtie H}\coloneqq(\mu_{A}\otimes\mu_{H})\circ (A\otimes\Psi\otimes H),\\
\varepsilon_{A\bowtie H}\coloneqq\varepsilon_{A\otimes H}, \;\;\; \delta_{A\bowtie H}\coloneqq\delta_{A\otimes H},\\
\lambda_{A\triangleright\hspace{-0.04 cm}\triangleleft H}\coloneqq \Psi\circ(\lambda_{H}\otimes\lambda_{A})\circ c_{A,H},
\end{gather*}
is a Hopf algebra called the bicrossed product of $A$ with $H$ (for details, see \cite[Theorem 7.2.2]{MAJ_Book}). Therefore, by the definition of $\mu_{A\bowtie H}$ and the unit property for $A$ and $H$, it is obtained that
\begin{equation}\label{psi_eta}
\Psi=\mu_{A\bowtie H}\circ(\eta_{A}\otimes H\otimes A\otimes\eta_{H}).
\end{equation}

\begin{remark}\label{particular_smash} 
Suppose that $(A,H,\varphi_{A},\phi_{H})$ is a matched pair of Hopf algebras such that $\phi_{H}=H\otimes\varepsilon_{A}$, i.e., $\phi_{H}$ is the trivial action of $A$ over $H$. This is equivalent to say that $(A,\varphi_{A})$ is a left $H$-module algebra-coalgebra such that the equality 
\begin{gather}\label{adcond_mp_phiHtriv}
(H\otimes\varphi_{A})\circ((c_{H,H}\circ\delta_{H})\otimes A)=(H\otimes\varphi_{A})\circ(\delta_{H}\otimes A)
\end{gather}
holds. 

Note that, in case $H$ is cocommutative, \eqref{adcond_mp_phiHtriv} is always satisfied. Following \cite{CCH}, this condition means that the pair $(A,\varphi_{A})$ is in the cocommutativity class of $H$ because, when we are in a symmetric context, \eqref{adcond_mp_phiHtriv} is equivalent to
\begin{equation}\label{cocom_class}
(\varphi_{A}\otimes H)\circ(H\otimes c_{H,A})\circ((c_{H,H}\circ\delta_{H})\otimes A)=(\varphi_{A}\otimes H)\circ(H\otimes c_{H,A})\circ(\delta_{H}\otimes A).
\end{equation}

When we  particularize the bicrossed product construction to this situation, we obtain the well-known smash product Hopf algebra of $A$ with $H$ which will be denoted by
\[A\sharp H=(A\otimes H,\eta_{A\otimes H},\mu_{A\sharp H},\varepsilon_{A\otimes H},\delta_{A\otimes H},\lambda_{A\sharp H}),\]
where
\begin{gather*}
\mu_{A\sharp H}\coloneqq (\mu_{A}\otimes\mu_{H})\circ (A\otimes((\varphi_{A}\otimes H)\circ(H\otimes c_{H,A})\circ(\delta_{H}\otimes A))\otimes H),\\
\lambda_{A\sharp H}\coloneqq (\varphi_{A}\otimes H)\circ(H\otimes c_{H,A})\circ(\delta_{H}\otimes A)\circ(\lambda_{H}\otimes\lambda_{A})\circ c_{A,H}.
\end{gather*}

In other words, the morphism $\Psi$ defined in \eqref{interw} becomes $\Psi=(\varphi_{A}\otimes H)\circ(H\otimes c_{H,A})\circ(\delta_{H}\otimes A)$ in this situation.

Moreover, if we also require that $\varphi_{A}=\varepsilon_{H}\otimes A$, then the usual tensor product Hopf algebra
\[A\otimes H=(A\otimes H,\eta_{A\otimes H},\mu_{A\otimes H},\varepsilon_{A\otimes H},\delta_{A\otimes H},\lambda_{A\otimes H})\]
is obtained, where $$\lambda_{A\otimes H}\coloneqq\lambda_{A}\otimes\lambda_{H},$$
i.e., $\Psi=c_{H,A}$ under such conditions.
\end{remark}

In what follows the structure of Hopf brace in {\sf C} is defined as well as giving some fundamental properties of these objects, which have been introduced by I. Angiono, C. Galindo and L. Vendramin in \cite{AGV} as the quantum version of skew braces (see \cite{GV}).
\begin{definition}\label{Hbrace}
Let $(H,\varepsilon_{H},\delta_{H})$ be a  coalgebra in {\sf C} and let's assume that $H$ admits two different algebra structures in {\sf C}: $(H,\eta_{H}^{1},\mu_{H}^{1})$ and $(H,\eta_{H}^{2},\mu_{H}^{2})$. We will say that a 9-tuple
\[(H,\eta_{H}^{1},\mu_{H}^{1},\eta_{H}^{2},\mu_{H}^{2},\varepsilon_{H},\delta_{H},\lambda_{H}^{1},\lambda_{H}^{2})\]
is a Hopf brace in {\sf C} if the following requirements hold:
\begin{itemize}
\item[(i)] $H_{1}=(H,\eta_{H}^{1},\mu_{H}^{1},\varepsilon_{H},\delta_{H},\lambda_{H}^{1})$ is a Hopf algebra in {\sf C}.
\item[(ii)]  $H_{2}=(H,\eta_{H}^{2},\mu_{H}^{2},\varepsilon_{H},\delta_{H},\lambda_{H}^{2})$ is a Hopf algebra in {\sf C}.
\item[(iii)] The following identity involving the products $\mu_{H}^{1}$ and $\mu_{H}^{2}$ holds:
\begin{gather}\label{compatHbrace}
\mu_{H}^{2}\circ(H\otimes\mu_{H}^{1})=\mu_{H}^{1}\circ(\mu_{H}^{2}\otimes\Gamma_{H_{1}})\circ(H\otimes c_{H,H}\otimes H)\circ(\delta_{H}\otimes H\otimes H),
\end{gather}
where 
\begin{equation}\label{def_GammaH1}
\Gamma_{H_{1}}\coloneqq \mu_{H}^{1}\circ(\lambda_{H}^{1}\otimes\mu_{H}^{2})\circ(\delta_{H}\otimes H).
\end{equation}
\end{itemize}

Following the notation introduced in \cite{RGON}, we will denote Hopf braces by $\mathbb{H}=(H_{1},H_{2})$ or, when there is no confusion between the Hopf algebras involved, only by $\mathbb{H}$. A Hopf brace $\mathbb{H}$ is said to be cocommutative if the underlying coalgebra structure, $(H,\varepsilon_{H},\delta_{H})$, is cocommutative, that is to say, if $c_{H,H}\circ\delta_{H}=\delta_{H}$.
\end{definition}

Given a Hopf brace $\mathbb{H}=(H_{1},H_{2})$, note that $$\eta_{H}^{1}=\eta_{H}^{2}$$ (see \cite[Remark 1.3]{AGV}). Therefore, from now on we will denote both units by $\eta_{H}$. Moreover, $(H_{1},\Gamma_{H_{1}})$ is a left $H_{2}$-module algebra and 
\begin{gather}\label{mu2-exp}
\mu_{H}^{2}=\mu_{H}^{1}\circ(H\otimes \Gamma_{H_{1}})\circ(\delta_{H}\otimes H),
\end{gather}
as was proved in \cite[Lemma 1.8]{AGV} and \cite[Remark 1.9]{AGV}, respectively. Note that \eqref{mu2-exp} holds without the need for \eqref{compatHbrace} to be satisfied, and then the following result sets a new equivalent characterization of Hopf braces.
\begin{theorem}\label{alt_Hbracecompat}
Let $H_{1}=(H,\eta_{H}^{1},\mu_{H}^{1},\varepsilon_{H},\delta_{H},\lambda_{H}^{1})$ and $H_{2}=(H,\eta_{H}^{2},\mu_{H}^{2},\varepsilon_{H},\delta_{H},\lambda_{H}^{2})$ be Hopf algebras in {\sf C}. Consider the morphism $\Gamma_{H_{1}}$ defined as in \eqref{def_GammaH1}. Then, the following statements are equivalent:
\begin{itemize}
\item[(i)] $\mathbb{H}=(H_{1},H_{2})$ is a Hopf brace in {\sf C},
\item[(ii)] $\Gamma_{H_{1}}\circ(H\otimes\mu_{H}^{1})=\mu_{H}^{1}\circ(\Gamma_{H_{1}}\otimes \Gamma_{H_{1}})\circ(H\otimes c_{H,H}\otimes H)\circ (\delta_{H}\otimes H\otimes H)$.
\end{itemize}
\end{theorem}
\begin{proof}
If $\mathbb{H}=(H_{1},H_{2})$ is assumed to be a Hopf brace, then $(H_{1},\Gamma_{H_{1}})$ is a left $H_{2}$-module algebra, what implies that (ii) holds. 

Suppose now that (ii) holds. As a consequence, 
\begin{itemize}
\itemindent=-32pt 
\item[ ]$\hspace{0.38cm}\mu_{H}^{1}\circ(\mu_{H}^{2}\otimes\Gamma_{H_{1}})\circ(H\otimes c_{H,H}\otimes H)\circ(\delta_{H}\otimes H\otimes H)$
\item[] $= \mu_{H}^{1}\circ ((\mu_{H}^{1}\circ(H\otimes\Gamma_{H_{1}})\circ(\delta_{H}\otimes H))\otimes\Gamma_{H_{1}})\circ(H\otimes c_{H,H}\otimes H)\circ(\delta_{H}\otimes H\otimes H)$ {\footnotesize\textnormal{(by \eqref{mu2-exp})}}
\item[] $= \mu_{H}^{1}\circ(H\otimes (\mu_{H}^{1}\circ(\Gamma_{H_{1}}\otimes \Gamma_{H_{1}})\circ(H\otimes c_{H,H}\otimes H)\circ (\delta_{H}\otimes H\otimes H)))\circ(\delta_{H}\otimes H\otimes H)$ {\footnotesize\textnormal{(by coassociativity}}
\item[] {\footnotesize\textnormal{of $\delta_{H}$ and associativity of $\mu_{H}^{1}$)}}
\item[] $=\mu_{H}^{1}\circ(H\otimes\Gamma_{H_{1}})\circ(\delta_{H}\otimes\mu_{H}^{1})$ {\footnotesize\textnormal{(by (ii))}}
\item[] $=\mu_{H}^{2}\circ(H\otimes\mu_{H}^{1})$ {\footnotesize\textnormal{(by \eqref{mu2-exp})}},
\end{itemize}
i.e., \eqref{compatHbrace} holds, and then $\mathbb{H}=(H_{1},H_{2})$ is a Hopf brace.
\end{proof}

Note that condition (ii) of Theorem \ref{alt_Hbracecompat} is equivalent to the fact that $(H_{1},\Gamma_{H_{1}})$ is a left $H_{2}$-module algebra because the equality $\Gamma_{H_{1}}\circ(H\otimes\eta_{H})=\varepsilon_{H}\otimes\eta_{H}$ is always satisfied. Also in \cite[Lemma 2.5]{VRBAProj} it is proved that $\Gamma_{H_{1}}$ is a coalgebra morphism if $\mathbb{H}$ is cocommutative, that is to say, $(H_{1},\Gamma_{H_{1}})$ is a left $H_{2}$-module algebra-coalgebra in the cocommutative setting.

\begin{remark}
By naturality of $c$, coassociativity of $\delta_{H}$ and associativity of $\mu_{H}^{1}$, 
\begin{align*}
&\mu_{H}^{1}\circ(\mu_{H}^{2}\otimes\Gamma_{H_{1}})\circ(H\otimes c_{H,H}\otimes H)\circ(\delta_{H}\otimes H\otimes H)\\=&\mu_{H}^{1}\circ(\Gamma'_{H_{1}}\otimes\mu_{H}^{2})\circ(H\otimes c_{H,H}\otimes H)\circ(\delta_{H}\otimes H\otimes H),
\end{align*}
where $$\Gamma'_{H_{1}}\coloneqq \mu_{H}^{1}\circ(\mu_{H}^{2}\otimes\lambda_{H}^{1})\circ(H\otimes c_{H,H})\circ(\delta_{H}\otimes H).$$ 

Therefore, condition (iii) of Definition \ref{Hbrace} is equivalent to
\begin{gather}\label{compatHbrace2}
\mu_{H}^{2}\circ(H\otimes\mu_{H}^{1})=\mu_{H}^{1}\circ(\Gamma'_{H_{1}}\otimes\mu_{H}^{2})\circ(H\otimes c_{H,H}\otimes H)\circ(\delta_{H}\otimes H\otimes H).
\end{gather}

Moreover, by \cite[Theorem 6.4]{Brz1}, $(H_{1},\Gamma'_{H_{1}})$ is a left $H_{2}$-module algebra and also $\Gamma'_{H_{1}}$ is a coalgebra morphism when $\mathbb{H}$ is a cocommutative Hopf brace (see \cite[Lemma 2.6]{VRBAProj}).

\end{remark}
\begin{definition}
Let $\mathbb{H}=(H_{1},H_{2})$ and $\mathbb{B}=(B_{1},B_{2})$ be Hopf braces in {\sf C}. A morphism $f\colon H\rightarrow B$ in {\sf C} between the underlying objects is said to be a morphism of Hopf braces if $f\colon H_{1}\rightarrow B_{1}$ and $f\colon H_{2}\rightarrow B_{2}$ are Hopf algebra morphisms.
\end{definition}

With the previous morphisms, Hopf braces constitute a category that we will denote by ${\sf HBr}$. Moreover, if we restrict ourselves to the cocommutative ones, then the full subcategory of cocommutative Hopf braces is obtained, which will be denoted by ${\sf cocHBr}$.

\begin{example}\label{trivHBr}
Every Hopf algebra $H=(H,\eta_{H},\mu_{H},\varepsilon_{H},\delta_{H},\lambda_{H})$ in {\sf C} gives rise to a Hopf brace structure by considering $H_{1}=H=H_{2}$. This Hopf brace will be called the trivial Hopf brace and we will denote it by $\mathbb{H}_{triv}=(H,H).$
\end{example}

It is important to recall that, under cocommutativity hypothesys, $H$ is a right $H_{2}$-module coalgebra with the action given by
\begin{gather*}
\Phi_{H}\coloneqq \mu_{H}^{2}\circ((\lambda_{H}^{2}\circ\Gamma_{H_{1}})\otimes\mu_{H}^{2})\circ(H\otimes c_{H,H}\otimes H)\circ(\delta_{H}\otimes\delta_{H}),
\end{gather*}
result whose proof can be seen in \cite[Lemma 2.2]{AGV}.

In \cite{AGV}, I. Angiono, C. Galindo and L. Vendramin also emphasize the strong relation between matched pairs of cocommutative Hopf algebras and cocommutative Hopf braces. Fixed $A$ a cocommutative Hopf algebra in {\sf C}, this connection will be detailed in what follows: On the one hand, if $\mathbb{A}=(A_{1},A)$ is a Hopf brace, then the 4-tuple $(A,A,\Gamma_{A_{1}},\Phi_{A})$ is a matched pair of Hopf algebras (see \cite[Proposition 3.1]{AGV}). Conversely, in \cite[Proposition 3.2]{AGV} it is proved that if $(A,A,\varphi_{A},\phi_{A})$ is a matched pair of Hopf algebras satisfying that 
\begin{gather}\label{cond_equiv_MP_HBrace}
\mu_{A}=\mu_{A}\circ\Psi,
\end{gather}
then $\overline{\mathbb{A}}=(\overline{A},A)$
is a Hopf brace, being 
\[\overline{A}=(A,\eta_{A},\mu_{\overline{A}},\varepsilon_{A},\delta_{A},\lambda_{\overline{A}})\]
the Hopf algebra whose product and antipode are defined as follows:
\begin{gather*}
\mu_{\overline{A}}\coloneqq\mu_{A}\circ(A\otimes(\varphi_{A}\circ(\lambda_{A}\otimes A)))\circ(\delta_{A}\otimes A),\\
\lambda_{\overline{A}}\coloneqq \varphi_{A}\circ(A\otimes \lambda_{A})\circ\delta_{A}.
\end{gather*}

The previous correspondence can be interpreted in a functorial way. Let's define ${\sf MP}(A)$ the category of matched pairs over the Hopf algebra $A$, whose objects are matched pairs $(A,A,\varphi_{A},\phi_{A})$ satisfying \eqref{cond_equiv_MP_HBrace}, while morphisms are defined as follows: $f\colon (A,A,\varphi_{A},\phi_{A})\rightarrow (A,A,\varphi'_{A},\phi'_{A})$ is a morphism in ${\sf MP}(A)$ if $f\colon A\rightarrow A$ is a Hopf algebra morphism satisfying that 
\[f\circ\varphi_{A}=\varphi'_{A}\circ(f\otimes f),\quad f\circ\phi_{A}=\phi'_{A}\circ(f\otimes f).\]

With ${\sf cocHBr}(A)$ we will denote the full subcategory of {\sf cocHBr} whose objects are Hopf braces of the form $\mathbb{A}=(A_{1},A)$, i.e., the second Hopf algebra structure is the fixed cocommutative Hopf algebra $A$. Therefore, we can define the following functors: 
\[F\colon {\sf cocHBr}(A)\longrightarrow {\sf MP}(A)\]
defined on objects by $F(\mathbb{A}=(A_{1},A))=(A,A,\Gamma_{A_{1}},\Phi_{A})$ and on morphisms by the identity, and
\[G\colon {\sf MP}(A)\longrightarrow {\sf cocHBr}(A)\]
defined on objects by $G((A,A,\varphi_{A},\phi_{A}))=\overline{\mathbb{A}}$ and on morphisms also by the identity. These functors give rise to a categorical isomorphism between the categories ${\sf cocHBr}(A)$ and ${\sf MP}(A)$ (see \cite[Theorem 3.3]{AGV}).

\section{Obtaining new Hopf braces from crossed products}
Motivated by the background mentioned in the last paragraphs of the previous section, it results natural to think about obtaining new Hopf braces using bicrossed products of Hopf algebras. Therefore, the problem we are going to study is the following: Let's suppose that $\mathbb{A}=(A_{1},A_{2})$ and $\mathbb{H}=(H_{1},H_{2})$ are Hopf braces in {\sf C} such that $(A_{1},H_{1},\varphi_{A}^{1},\phi_{H}^{1})$ and $(A_{2},H_{2},\varphi_{A}^{2},\phi_{H}^{2})$ are matched pairs of Hopf algebras. Then 
$$A_{1}\bowtie H_{1}=(A\otimes H,\eta_{A\otimes H},\mu_{A\bowtie H}^{1},\varepsilon_{A\otimes H},\delta_{A\otimes H},\lambda_{A\bowtie H}^{1}),$$
$$A_{2}\bowtie H_{2}=(A\otimes H,\eta_{A\otimes H},\mu_{A\bowtie H}^{2},\varepsilon_{A\otimes H},\delta_{A\otimes H},\lambda_{A\bowtie H}^{2})$$
are Hopf algebras in {\sf C} with the same underlying coalgebra structure, where 
\begin{gather*}
\mu_{A\bowtie H}^{i}=(\mu_{A}^i\otimes\mu_{H}^i)\circ (A\otimes\Psi^{i}\otimes H),\\
\lambda_{A\bowtie H}^{i}=\Psi^{i}\circ(\lambda_{H}^{i}\otimes\lambda_{A}^{i})\circ c_{A,H},
\end{gather*}
being $\Psi^{i}=(\varphi_{A}^{i}\otimes\phi_{H}^{i})\circ(H\otimes c_{H,A}\otimes A)\circ(\delta_{H}\otimes\delta_{A})$ for all $i=1,2$. So, the question we want to address is the following:
\begin{equation}\label{prob}\tag{P}
\textnormal{Under what conditions $(A_{1}\bowtie H_{1},A_{2}\bowtie H_{2})$ gives rise to a new Hopf brace in {\sf C}?}
\end{equation} 

Despite giving a solution to this problem is not an easy task, A. Agore in \cite{AGORE} has been able to tackle it in certain particular cases that will be mentioned hereunder. This first result can be consulted in \cite[Theorem 2.1]{AGORE}.
\begin{theorem}\label{agtheorem1}
Let $\mathbb{A}=(A_{1},A_{2})$ be a Hopf brace and $H$ a commutative and cocommutative Hopf algebra in {\sf C}. If the conditions
\begin{itemize}
\item[(i)] $(A_{1},\varphi_{A})$ is a left $H$-module algebra,
\item[(ii)] $(A_{2},H,\varphi_{A},\phi_{H})$ is a matched pair of Hopf algebras,
\item[(iii)] $\phi_{H}\circ (H\otimes\mu_{A}^{1})=\mu_{H}\circ (\phi_{H}\otimes (\mu_{H}\circ (\lambda_{H}\otimes\phi_{H})\circ (\delta_{H}\otimes A)))\circ (H\otimes c_{H,A}\otimes A)\circ(\delta_{H}\otimes A\otimes A)$
\end{itemize}
hold, then 
$(A_{1}\otimes H,A_{2}\bowtie H)$
is a Hopf brace in {\sf C}.
\end{theorem}
The second result proved by A. Agore in \cite[Theorem 2.5]{AGORE} is the following:
\begin{theorem}\label{agtheorem2}
Let $A$ be a Hopf algebra and $\mathbb{H}=(H_{1},H_{2})$ a cocommutative Hopf brace in {\sf C}. If the conditions
\begin{itemize}
\item[(i)] $(A,H_{1},\varphi_{A}^{1},\phi_{H})$ is a matched pair of Hopf algebras,
\item[(ii)] $(A,\varphi_{A}^{2})$ is a left $H_{2}$-module algebra-coalgebra,
\item[(iii)] $\varphi_{A}^{2}\circ(H\otimes\varphi_{A}^{1})=\varphi_{A}^{1}\circ (\Gamma'_{H_{1}}\otimes \varphi_{A}^{2})\circ (H\otimes c_{H,H}\otimes A)\circ(\delta_{H}\otimes H\otimes A)$,
\item[(iv)] $\mu_{H}^{2}\circ(H\otimes\phi_{H})=\mu_{H}^{1}\circ(\phi_{H}\otimes H)\circ (\Gamma'_{H_{1}}\otimes\Psi^{2})\circ (H\otimes c_{H,H}\otimes A)\circ (\delta_{H}\otimes H\otimes A)$
\end{itemize}
hold, then
$(A\bowtie H_{1},A\sharp H_{2})$
is a Hopf brace in {\sf C}.
\end{theorem}
\begin{remark}
Considering the definition of module over a Hopf brace introduced by R. González in \cite[Definition 2.10]{RGON}, hypothesis (iii) of previous theorem is equivalent to the fact that $(A,\varphi_{A}^{1},\varphi_{A}^{2})$ is a left $\mathbb{H}$-module.
\end{remark}

So, in the previous theorems, A. Agore combines a Hopf algebra with a Hopf brace using crossed products in order to obtain new Hopf brace structures assuming some additional hypothesis, commutativity or cocommutativity always included. Although we did not manage to obtain a complete solution for \eqref{prob}, the main results of this article (Theorems \ref{mainth2} and \ref{mainth}) go one step further than A. Agore's results because they consist of the first constructions of new Hopf braces mixing two different Hopf braces using crossed products.
\begin{theorem}\label{mainth2}
Let $\mathbb{A}=(A_{1},A_{2})$ be an object in {\sf HBr} and $\mathbb{H}=(H_{1},H_{2})$ be an object in {\sf cocHBr} such that $(A_{2},H_{2},\varphi_{A},\phi_{H})$ is a matched pair of Hopf algebras. Under these conditions, the pair $$(A_{1}\otimes H_{1},A_{2}\bowtie H_{2})$$ is a Hopf brace in {\sf C} if and only if the following conditions hold:
\begin{itemize}
\item[(E1)] $\mu_{A}^{1}\circ (\varphi_{A}\otimes\varphi_{A})\circ (H\otimes c_{H,A}\otimes A)\circ(\delta_{H}\otimes A\otimes A)=\varphi_{A}\circ(H\otimes\mu_{A}^{1})$,
\item[(E2)] $\hspace{0.38cm}\mu_{H}^{1}\circ (\mu_{H}^{2}\otimes \Omega)\circ(H\otimes c_{H,H}\otimes A\otimes H)\circ (((\phi_{H}\otimes H)\circ(H\otimes c_{H,A})\circ (\delta_{H}\otimes A))\otimes H\otimes A\otimes H)\\=\mu_{H}^{2}\circ((\phi_{H}\circ(H\otimes\mu_{A}^{1}))\otimes\mu_{H}^{1})\circ(H\otimes A\otimes c_{H,A}\otimes H)$,
\end{itemize}
where 
\begin{equation}\label{Omega_def}\Omega\coloneqq \mu_{H}^{1}\circ (\lambda_{H}^{1}\otimes (\mu_{H}^{2}\circ (\phi_{H}\otimes H)))\circ (\delta_{H}\otimes A\otimes H).\end{equation}
\end{theorem}
\begin{proof} At first we will prove some necessary equalities for the proof, which are the followings:
\begin{gather}
\label{e1} (\varepsilon_{A}\otimes H)\circ \Psi=\phi_{H},\\
\label{e2}(A\otimes\varepsilon_{H})\circ\Psi=\varphi_{A},\\
\label{e3} \Gamma_{A_{1}\otimes H_{1}}=((\Gamma_{A_{1}}\circ(A\otimes\varphi_{A}))\otimes\Omega)\circ(A\otimes\delta_{H\otimes A}\otimes H),\\
\label{e4} (\varepsilon_{A}\otimes H)\circ \Gamma_{A_{1}\otimes H_{1}}\circ (\eta_{A}\otimes H\otimes A\otimes H)=\Omega,\\
\label{e5} (A\otimes\varepsilon_{H})\circ \Gamma_{A_{1}\otimes H_{1}}\circ(A\otimes H\otimes A\otimes\eta_{H})=\Gamma_{A_{1}}\circ(A\otimes\varphi_{A}).
\end{gather}

Equalities \eqref{e1} and \eqref{e2} are straightforward to prove using the condition of coalgebra morphism for $\varphi_{A}$ and $\phi_{H}$, respectively, naturality of $c$ and the counit properties. Moreover, \eqref{e3} follows by
\begin{itemize}
\itemindent=-32pt 
\item[ ]$\hspace{0.38cm} \Gamma_{A_{1}\otimes H_{1}}$
\item[]$=(\Gamma_{A_{1}}\otimes (\mu_{H}^{1}\circ(\lambda_{H}^{1}\otimes \mu_{H}^{2})))\circ(A\otimes c_{H,A}\otimes H\otimes H)\circ(A\otimes H\otimes\Psi\otimes H)\circ(A\otimes \delta_{H}\otimes A\otimes H)$ {\footnotesize (by} 
\item[] {\footnotesize naturality of $c$ and the symmetric character of {\sf C})}
\item[]$=((\Gamma_{A_{1}}\circ(A\otimes\varphi_{A}))\otimes(\mu_{H}^{1}\circ (\lambda_{H}^{1}\otimes\mu_{H}^{2})\circ (H\otimes\phi_{H}\otimes H)))\circ(A\otimes((H\otimes c_{H,A}\otimes H)\circ((c_{H,H}\circ\delta_{H})$
\item[]$\otimes c_{H,A}))\otimes A\otimes H)\circ(A\otimes\delta_{H}\otimes\delta_{A}\otimes H)$ {\footnotesize\textnormal{(by naturality of $c$ and coassociativity of $\delta_{H}$)}}
\item[]$=((\Gamma_{A_{1}}\circ(A\otimes\varphi_{A}))\otimes(\mu_{H}^{1}\circ (\lambda_{H}^{1}\otimes\mu_{H}^{2})\circ (H\otimes\phi_{H}\otimes H)))\circ(A\otimes H\otimes ((c_{H,A}\otimes H)\circ(H\otimes c_{H,A})\circ$
\item[] $(\delta_{H}\otimes A))\otimes A\otimes H)\circ(A\otimes\delta_{H}\otimes\delta_{A}\otimes H)$ {\footnotesize\textnormal{(by cocommutativity and coassociativity of $\delta_{H}$)}}
\item[] $=((\Gamma_{A_{1}}\circ(A\otimes\varphi_{A}))\otimes\Omega)\circ(A\otimes\delta_{H\otimes A}\otimes H)$ {\footnotesize\textnormal{(by naturality of $c$)}.}
\end{itemize}

As a consequence of \eqref{e3}, equalities \eqref{e4} and \eqref{e5} hold. Indeed, on the one hand, \eqref{e4} follows from \eqref{e3} by using that $\varepsilon_{A}\circ \Gamma_{A_{1}}=\varepsilon_{A}\otimes\varepsilon_{A}$, the condition of coalgebra morphism for $\varphi_{A}$, the counit property for $A$ and $H$ and the equality $\varepsilon_{A}\circ\eta_{A}=id_{K}$. On the other hand, \eqref{e5} is deduced from \eqref{e3} by applying the condition of coalgebra morphism for $\mu_{H}^{1},$ $\mu_{H}^{2}$ and $\phi_{H}$, \eqref{u-antip2}, the counit property for $H$ and $A$ and the equality $\varepsilon_{H}\circ\eta_{H}=id_{K}$.

Once we have justified equations \eqref{e1}-\eqref{e5}, let's start the proof of the equivalence by assuming that the pair $(A_{1}\otimes H_{1},A_{2}\bowtie H_{2})$ is a Hopf brace in {\sf C}, what implies that the following equality holds:
\begin{align}\label{compatbraceparticular}
&\mu_{A\bowtie H}^{2}\circ (A\otimes H\otimes\mu_{A\otimes H}^{1})\\\nonumber=&\mu_{A\otimes H}^{1}\circ(\mu_{A\bowtie H}^{2}\otimes \Gamma_{A_{1}\otimes H_{1}})\circ(A\otimes H\otimes c_{A\otimes H,A\otimes H}\otimes A\otimes H)\circ(\delta_{A\otimes H}\otimes A\otimes H\otimes A\otimes H). 
\end{align}

On the one side, composing on the right of \eqref{compatbraceparticular} with $\eta_{A}\otimes H\otimes A\otimes \eta_{H}\otimes A\otimes\eta_{H}$, and on the left with $A\otimes\varepsilon_{H}$, from the left hand side it follows that
\begin{itemize}
\itemindent=-32pt 
\item[ ]$\hspace{0.38cm} (A\otimes\varepsilon_{H})\circ \mu_{A\bowtie H}^{2}\circ (A\otimes H\otimes\mu_{A\otimes H}^{1})\circ (\eta_{A}\otimes H\otimes A\otimes \eta_{H}\otimes A\otimes\eta_{H})$
\item[]$=(A\otimes\varepsilon_{H})\circ\Psi\circ(H\otimes\mu_{A}^{1})$ {\footnotesize\textnormal{(by the condition of coalgebra morphism for $\mu_{H}^{1}$ and $\mu_{H}^{2}$, the equality $\varepsilon_{H}\circ\eta_{H}=id_{K}$}} 
\item[] {\footnotesize\textnormal{and the unit property for $A_{2}$)}}
\item[] $=\varphi_{A}\circ(H\otimes\mu_{A}^{1})$ {\footnotesize\textnormal{(by \eqref{e2})},}
\end{itemize}
whereas, from the right hand side, it is obtained that
\begin{itemize}
\itemindent=-32pt 
\item[ ]$\hspace{0.38cm}(A\otimes\varepsilon_{H})\circ\mu_{A\otimes H}^{1}\circ(\mu_{A\bowtie H}^{2}\otimes \Gamma_{A_{1}\otimes H_{1}})\circ(A\otimes H\otimes c_{A\otimes H,A\otimes H}\otimes A\otimes H)\circ((\delta_{A\otimes H}\circ(\eta_{A}\otimes H))\otimes $
\item[] $A\otimes \eta_{H}\otimes A\otimes \eta_{H})$
\item[] $=\mu_{A}^{1}\circ (((A\otimes\varepsilon_{H})\circ\Psi)\otimes((A\otimes\varepsilon_{H})\circ \Gamma_{A_{1}\otimes H_{1}}\circ (\eta_{A}\otimes H\otimes A\otimes\eta_{H})))\circ (H\otimes c_{H,A}\otimes A)\circ(\delta_{H}\otimes A\otimes A)$ 
\item []{\footnotesize\textnormal{(by the condition of coalgebra morphism for $\mu_{H}^{1}$ and $\mu_{H}^{2}$, the condition of coalgebra morphism for $\eta_{A}$, the naturality}}
\item[]{\footnotesize\textnormal{of $c$, the unit property for $A_{2}$ and the equality $\varepsilon_{H}\circ\eta_{H}=id_{K}$)}}
\item[] $=\mu_{A}^{1}\circ (\varphi_{A}\otimes (\Gamma_{A_{1}}\circ(\eta_{A}\otimes\varphi_{A})))\circ(H\otimes c_{H,A}\otimes A)\circ(\delta_{H}\otimes A\otimes A)$ {\footnotesize\textnormal{(by \eqref{e2} and \eqref{e5})}}
\item[] $=\mu_{A}^{1}\circ(\varphi_{A}\otimes\varphi_{A})\circ(H\otimes c_{H,A}\otimes A)\circ(\delta_{H}\otimes A\otimes A)$ {\footnotesize\textnormal{(by module axioms for $(A_{1},\Gamma_{A_{1}})$)}.}
\end{itemize}

As a consequence, we deduce that $\mu_{A}^{1}$ is a morphism of left $H$-modules, i.e., (E1) holds.

On the other side, composing on the right of \eqref{compatbraceparticular} with $\eta_{A}\otimes H\otimes A\otimes H\otimes A\otimes H$, and on the left with $\varepsilon_{A}\otimes H$, from the left hand side it follows that
\begin{itemize}
\itemindent=-32pt 
\item[ ]$\hspace{0.38cm}(\varepsilon_{A}\otimes H)\circ\mu_{A\bowtie H}^{2}\circ (A\otimes H\otimes\mu_{A\otimes H}^{1})\circ(\eta_{A}\otimes H\otimes A\otimes H\otimes A\otimes H)$
\item[] $=\mu_{H}^{2}\circ(((\varepsilon_{A}\otimes H)\circ\Psi)\otimes H)\circ(H\otimes((\mu_{A}^{1}\otimes\mu_{H}^{1})\circ (A\otimes c_{H,A}\otimes H)))$ {\footnotesize\textnormal{(by the condition of coalgebra}}
\item[] {\footnotesize\textnormal{morphism for $\mu_{A}^{2}$ and $\varepsilon_{A}\circ\eta_{A}=id_{K}$)}}
\item[] $=\mu_{H}^{2}\circ((\phi_{H}\circ(H\otimes\mu_{A}^{1}))\otimes\mu_{H}^{1})\circ(H\otimes A\otimes c_{H,A}\otimes H)$ {\footnotesize\textnormal{(by \eqref{e1})},}
\end{itemize}
while, from the right hand side, it is obtained that
\begin{itemize}
\itemindent=-32pt 
\item[ ]$\hspace{0.38cm}(\varepsilon_{A}\otimes H)\circ\mu_{A\otimes H}^{1}\circ(\mu_{A\bowtie H}^{2}\otimes \Gamma_{A_{1}\otimes H_{1}})\circ(A\otimes H\otimes c_{A\otimes H,A\otimes H}\otimes A\otimes H)\circ((\delta_{A\otimes H}\circ(\eta_{A}\otimes H))\otimes $
\item[] $A\otimes H\otimes A\otimes H)$
\item[ ] $=\mu_{H}^{1}\circ (\mu_{H}^{2}\otimes((\varepsilon_{A}\otimes H)\circ\Gamma_{A_{1}\otimes H_{1}}\circ(\eta_{A}\otimes H\otimes A\otimes H)))\circ (((\varepsilon_{A}\otimes H)\circ\Psi)\otimes c_{H,H}\otimes A\otimes H)\circ $
\item[] $(((H\otimes c_{H,A})\circ(\delta_{H}\otimes A))\otimes H\otimes A\otimes H)$ {\footnotesize\textnormal{(by the condition of coalgebra morphism for $\mu_{A}^{1}$ and $\mu_{A}^{2}$, the naturality}}
\item[] {\footnotesize\textnormal{of $c$ and the counit property for $A$)}}
\item[] $=\mu_{H}^{1}\circ (\mu_{H}^{2}\otimes\Omega)\circ (H\otimes c_{H,H}\otimes A\otimes H)\circ(((\phi_{H}\otimes H)\circ (H\otimes c_{H,A})\circ (\delta_{H}\otimes A))\otimes H\otimes A\otimes H)$ {\footnotesize\textnormal{(by}}
\item[] {\footnotesize\textnormal{\eqref{e1} and \eqref{e4})},}
\end{itemize}
what implies that (E2) also holds. Then, the first implication is done.

Now, let's suppose that (E1) and (E2) hold. Therefore,
\begin{itemize}
\itemindent=-32pt 
\item[ ]$\hspace{0.38cm}\mu_{A\otimes H}^{1}\circ(\mu_{A\bowtie H}^{2}\otimes \Gamma_{A_{1}\otimes H_{1}})\circ(A\otimes H\otimes c_{A\otimes H,A\otimes H}\otimes A\otimes H)\circ(\delta_{A\otimes H}\otimes A\otimes H\otimes A\otimes H)$
\item[] $=(\mu_{A}^{1}\otimes\mu_{H}^{1})\circ(A\otimes c_{H,A}\otimes H)\circ(((\mu_{A}^{2}\otimes\mu_{H}^{2})\circ(A\otimes\Psi\otimes H))\otimes (((\Gamma_{A_{1}}\circ(A\otimes\varphi_{A}))\otimes\Omega)\circ(A\otimes$
\item[] $\delta_{H\otimes A}\otimes H)))\circ(A\otimes H\otimes c_{A\otimes H,A\otimes H}\otimes A\otimes H)\circ(\delta_{A\otimes H}\otimes A\otimes H\otimes A\otimes H)$ {\footnotesize\textnormal{(by \eqref{e3})}}
\item[] $=((\mu_{A}^{1}\circ(\mu_{A}^{2}\otimes\Gamma_{A_{1}})\circ(A\otimes c_{A,A}\otimes\varphi_{A})\circ(\delta_{A}\otimes A\otimes H\otimes A))\otimes(\mu_{H}^{1}\circ (\mu_{H}^{2}\otimes\Omega)))\circ(A\otimes((A\otimes c_{H\otimes H,H\otimes A}\otimes$
\item[]$ H\otimes A)\circ(A\otimes H\otimes H\otimes\delta_{H\otimes A})\circ(\Psi\otimes c_{H,H}\otimes A)\circ(H\otimes c_{H,A}\otimes H\otimes A)\circ(\delta_{H}\otimes A\otimes H\otimes A))\otimes H)$ 
\item[] {\footnotesize\textnormal{(by naturality of $c$ and symmetrical character of {\sf C})}}
\item[] $=((\mu_{A}^{2}\circ(A\otimes\mu_{A}^{1})\circ(A\otimes A\otimes \varphi_{A}))\otimes(\mu_{H}^{1}\circ (\mu_{H}^{2}\otimes\Omega)))\circ(A\otimes((A\otimes c_{H\otimes H,H\otimes A}\otimes H\otimes A)\circ(A\otimes H\otimes$
\item[]$ H\otimes\delta_{H\otimes A})\circ(\Psi\otimes c_{H,H}\otimes A)\circ(H\otimes c_{H,A}\otimes H\otimes A)\circ(\delta_{H}\otimes A\otimes H\otimes A))\otimes H)$ {\footnotesize\textnormal{(by \eqref{compatHbrace} for $\mathbb{A}$)}}
\item[] $=((\mu_{A}^{2}\circ(A\otimes(\mu_{A}^{1}\circ(\varphi_{A}\otimes\varphi_{A})\circ(H\otimes c_{H,A}\otimes A))))\otimes(\mu_{H}^{1}\circ ((\mu_{H}^{2}\circ(\phi_{H}\otimes H))\otimes\Omega)\circ(H\otimes A\otimes c_{H,H}\otimes A\otimes$
\item[] $ H)))\circ(A\otimes H\otimes H\otimes((A\otimes A\otimes H\otimes c_{H,A}\otimes H)\circ(c_{H\otimes H,A\otimes A}\otimes A\otimes H)\circ(H\otimes H\otimes A\otimes((c_{A,A}\otimes H)\circ$
\item[]$(A\otimes c_{H,A}))))\otimes A\otimes H)\circ(A\otimes((H\otimes(c_{H,H}\circ\delta_{H})\otimes H)\circ(H\otimes\delta_{H})\circ\delta_{H})\otimes\delta_{A}\otimes H\otimes\delta_{A}\otimes H)$ {\footnotesize\textnormal{(by}}\item[] {\footnotesize\textnormal{naturality of $c$, the symmetric character of {\sf C} and coassociativity of $\delta_{H}$)}}
\item[]$=((\mu_{A}^{2}\circ(A\otimes(\mu_{A}^{1}\circ(\varphi_{A}\otimes\varphi_{A})\circ(H\otimes c_{H,A}\otimes A)\circ(\delta_{H}\otimes A\otimes A))))\otimes(\mu_{H}^{1}\circ (\mu_{H}^{2}\otimes \Omega)\circ(H\otimes c_{H,H}\otimes A\otimes$
\item[]$ H)\circ (((\phi_{H}\otimes H)\circ(H\otimes c_{H,A})\circ (\delta_{H}\otimes A))\otimes H\otimes A\otimes H)))\circ(A\otimes H\otimes A\otimes((c_{H,A}\otimes A\otimes H\otimes A)\circ(H\otimes $
\item[] $c_{A,A}\otimes H\otimes A)\circ(H\otimes A\otimes c_{H,A}\otimes A)\circ(H\otimes A\otimes H\otimes\delta_{A}))\otimes H)\circ(A\otimes\delta_{H\otimes A}\otimes H\otimes A\otimes H)$ {\footnotesize\textnormal{(by}}
\item[] {\footnotesize\textnormal{cocommutativity and coassociativity of $\delta_{H}$ and naturality of $c$)}}
\item[] $=((\mu_{A}^{2}\circ(A\otimes(\varphi_{A}\circ(H\otimes\mu_{A}^{1}))))\otimes(\mu_{H}^{2}\circ((\phi_{H}\circ(H\otimes\mu_{A}^{1}))\otimes\mu_{H}^{1})\circ(H\otimes A\otimes c_{H,A}\otimes H)))\circ(A\otimes $
\item[] $H\otimes A\otimes((c_{H,A}\otimes A\otimes H\otimes A)\circ(H\otimes c_{A,A}\otimes H\otimes A)\circ(H\otimes A\otimes c_{H,A}\otimes A)\circ(H\otimes A\otimes H\otimes\delta_{A}))\otimes$
\item[] $ H)\circ(A\otimes\delta_{H\otimes A}\otimes H\otimes A\otimes H)$ {\footnotesize\textnormal{(by (E1) and (E2))}}
\item[] $=(\mu_{A}^{2}\otimes\mu_{H}^{2})\circ (A\otimes((\varphi_{A}\otimes\phi_{H})\circ(H\otimes c_{H,A}\otimes A)\circ(\delta_{H}\otimes((\mu_{A}^{1}\otimes\mu_{A}^{1})\circ(A\otimes c_{A,A}\otimes A)\circ(\delta_{A}\otimes\delta_{A}))))\otimes $
\item[] $\mu_{H}^{1})\circ (A\otimes H\otimes A\otimes c_{H,A}\otimes H)$ {\footnotesize\textnormal{(by naturality of $c$)}}
\item[] $=\mu_{A\bowtie H}^{2}\circ(A\otimes H\otimes\mu_{A\otimes H}^{1})$ {\footnotesize\textnormal{(by the condition of coalgebra morphism for $\mu_{A}^{1}$)},}
\end{itemize}
what implies that $(A_{1}\otimes H_{1},A_{2}\bowtie H_{2})\in{\sf HBr}$.
\end{proof}
Note that condition (E1) of Theorem \ref{mainth2} is equivalent to the fact that $(A_{1},\varphi_{A})$ is a left $H_{2}$-module algebra because the equality $\varphi_{A}\circ(H\otimes\eta_{A})=\varepsilon_{H}\otimes\eta_{A}$ automatically holds by (ii) of Definition \ref{MPdef}.

As a particular situation of the previous theorem, A. Agore's Theorem \ref{agtheorem1} is deduced.
\begin{corollary}
Let $\mathbb{A}=(A_{1},A_{2})$ be a Hopf brace and $H$ a commutative and cocommutative Hopf algebra in {\sf C}. If the conditions
\begin{itemize}
\item[(i)] $(A_{1},\varphi_{A})$ is a left $H$-module algebra,
\item[(ii)] $(A_{2},H,\varphi_{A},\phi_{H})$ is a matched pair of Hopf algebras,
\item[(iii)] $\phi_{H}\circ (H\otimes\mu_{A}^{1})=\mu_{H}\circ (\phi_{H}\otimes (\mu_{H}\circ (\lambda_{H}\otimes\phi_{H})\circ (\delta_{H}\otimes A)))\circ (H\otimes c_{H,A}\otimes A)\circ(\delta_{H}\otimes A\otimes A)$
\end{itemize}
hold, then 
$(A_{1}\otimes H,A_{2}\bowtie H)$
is a Hopf brace in {\sf C}.
\end{corollary}
\begin{proof}
If we consider $\mathbb{H}=\mathbb{H}_{triv}=(H,H)$ the trivial Hopf brace built from the commutative and cocommutative Hopf algebra $H$, then we are in the conditions of Theorem \ref{mainth2}. Thus, in order to prove that $(A_{1}\otimes H, A_{2}\bowtie H)$ is a Hopf brace we only need to compute that conditions (E1) and (E2) hold in this situation. Note that (E1) is equivalent to the condition (i), while (E2) is tested as follows:
\begin{itemize}
\itemindent=-32pt 
\item[ ]$\hspace{0.38cm}\mu_{H}\circ (\mu_{H}\otimes \Omega)\circ(H\otimes c_{H,H}\otimes A\otimes H)\circ (((\phi_{H}\otimes H)\circ(H\otimes c_{H,A})\circ (\delta_{H}\otimes A))\otimes H\otimes A\otimes H)$
\item[] $=\mu_{H}\circ((\mu_{H}\circ (H\otimes (\mu_{H}\circ c_{H,H}\circ (H\otimes(\mu_{H}\circ (\lambda_{H}\otimes\phi_{H})\circ(\delta_{H}\otimes A)))\circ(c_{H,H}\otimes A))))\otimes H)\circ(((\phi_{H}\otimes$ 
\item[] $H)\circ(H\otimes c_{H,A})\circ (\delta_{H}\otimes A))\otimes H\otimes A\otimes H)$ {\footnotesize\textnormal{(by associativity and commutativity of $\mu_{H}$)}}
\item[] $=\mu_{H}\circ ((\mu_{H}\circ (\phi_{H}\otimes (\mu_{H}\circ (\lambda_{H}\otimes\phi_{H})\circ(\delta_{H}\otimes A)))\circ (H\otimes c_{H,A}\otimes A)\circ(\delta_{H}\otimes A\otimes A))\otimes\mu_{H})\circ(H\otimes A\otimes c_{H,A}\otimes H)$ 
\item[] {\footnotesize\textnormal{(by naturality of $c$, symmetrical character of {\sf C} and associativity of $\mu_{H}$)}}
\item[] $=\mu_{H}\circ ((\phi_{H}\circ (H\otimes\mu_{A}^{1}))\otimes\mu_{H})\circ(H\otimes A\otimes c_{H,A}\otimes H)$ {\footnotesize\textnormal{(by (iii)).}}\qedhere
\end{itemize}
\end{proof}
\begin{remark}
In the proof of the previous corollary we have seen that hypothesis (iii) implies condition (E2) of Theorem \ref{mainth2}. However, it is in fact an equivalence, that is to say, (E2) also implies (iii). To prove this fact it is enough to compose on the right of equality (E2) with $H\otimes A\otimes\eta_{H}\otimes A\otimes\eta_{H}$.

\begin{corollary}
Let $A$ be a Hopf algebra and $H$ a cocommutative Hopf algebra in {\sf C}. If $(A,H,\varphi_{A},\phi_{H})$ is a matched pair of Hopf algebras, then $(A\otimes H,A\bowtie H)$ is a Hopf brace in {\sf C} if and only if the equalities
\begin{align}
\label{mod_alg_z1}&\mu_{A}\circ(\varphi_{A}\otimes\varphi_{A})\circ(H\otimes c_{H,A}\otimes A)\circ(\delta_{H}\otimes A\otimes A)=\varphi_{A}\circ(H\otimes\mu_{A}),\\
\label{z2}&\mu_{H}\circ(\mu_{H}\otimes \overline{\Omega})\circ(H\otimes c_{H,H}\otimes A)\circ(((\phi_{H}\otimes H)\circ(H\otimes c_{H,A})\circ(\delta_{H}\otimes A))\otimes H\otimes A)\\\nonumber=&\mu_{H}\circ((\phi_{H}\circ(H\otimes\mu_{A}))\otimes H)\circ(H\otimes A\otimes c_{H,A})
\end{align}
hold, where 
\begin{equation}\label{Omegabar_def}\overline{\Omega}\coloneqq\mu_{H}\circ(\lambda_{H}\otimes\phi_{H})\circ(\delta_{H}\otimes A).\end{equation}
\end{corollary}
\begin{proof}
It is enough to apply Theorem \ref{mainth2} by considering that $\mathbb{H}=\mathbb{H}_{triv}=(H,H)$, the trivial cocommutative Hopf brace built from the cocommutative Hopf algebra $H$, and $\mathbb{A}=\mathbb{A}_{triv}=(A,A)$, the trivial Hopf brace obtained from the Hopf algebra $A$. In this situation, it results clear that (E1) is equivalent to \eqref{mod_alg_z1}, while \eqref{z2} is equivalent to (E2). Indeed, note that the morphism $\Omega$ defined in \eqref{Omega_def} satisfies that 
\begin{equation}\label{z3}
\Omega=\mu_{H}\circ(\overline{\Omega}\otimes H)
\end{equation}
by associativity of $\mu_{H}$. Then, on the one hand, if (E2) holds, \eqref{z2} is obtained by composing on the right with $H\otimes A\otimes H\otimes A\otimes \eta_{H}$ and, on the other hand, if \eqref{z2} holds, then
\begin{itemize}
\itemindent=-32pt 
\item[ ]$\hspace{0.38cm}\mu_{H}\circ (\mu_{H}\otimes \Omega)\circ(H\otimes c_{H,H}\otimes A\otimes H)\circ (((\phi_{H}\otimes H)\circ(H\otimes c_{H,A})\circ (\delta_{H}\otimes A))\otimes H\otimes A\otimes H)$
\item[] $=\mu_{H}\circ((\mu_{H}\circ(\mu_{H}\otimes\overline{\Omega})\circ (H\otimes c_{H,H}\otimes A)\circ(((\phi_{H}\otimes H)\circ(H\otimes c_{H,A})\circ(\delta_{H}\otimes A))\otimes H\otimes A)) \otimes H)$ 
\item[] {\footnotesize\textnormal{(by \eqref{z3} and associativity of $\mu_{H}$)}}
\item[] $=\mu_{H}\circ((\mu_{H}\circ((\phi_{H}\circ(H\otimes\mu_{A}))\otimes H)\circ(H\otimes A\otimes c_{H,A}))\otimes H)$ {\footnotesize\textnormal{(by \eqref{z2})}}
\item[] $=\mu_{H}\circ(\phi_{H}\otimes H)\circ(H\otimes((\mu_{A}\otimes\mu_{H})\circ (A\otimes c_{H,A}\otimes H)))$ {\footnotesize\textnormal{(by associativity of $\mu_{H}$)}},
\end{itemize}
what implies that (E2) is satisfied.
\end{proof}
\begin{remark} Under the conditions of previous corollary, if it is assumed that $H$ is also commutative, then it is straightforward to prove that (\ref{z2}) is equivalent to 
\[\mu_{H}\circ (\phi_{H}\otimes \overline{\Omega})\circ(H\otimes c_{H,A}\otimes A)\circ(\delta_{H}\otimes A\otimes A)=\phi_{H}\circ(H\otimes\mu_{A}).\]
\end{remark}

Looking at the statement of the above corollary, a natural question one may ask is the following: Is it possible to reduce (or even eliminate) the cocommutativity hypothesis for $H$? In the following theorem we will see that it is enough to require that $(A,\varphi_{A})$ is in the cocommutativity class of $H$ (see Remark \ref{particular_smash}).
\begin{theorem}\label{class_cocom_bicros} 
Let $A$ and $H$ be Hopf algebras in {\sf C}. If $(A,H,\varphi_{A},\phi_{H})$ is a matched pair of Hopf algebras such that $(A,\varphi_{A})$ is in the cocommutativity class of $H$, i.e., equality \eqref{cocom_class} holds, then $(A\otimes H,A\bowtie H)$ is a Hopf brace in {\sf C} if and only if equations \eqref{mod_alg_z1} and \eqref{z2} hold, where $\overline{\Omega}$ is defined as in \eqref{Omegabar_def}.
\end{theorem}
\begin{proof}
Let's start by computing $\Gamma_{A\otimes H}$:
\begin{itemize}
\itemindent=-32pt 
\item[ ]$\hspace{0.38cm} \Gamma_{A\otimes H}$
\item[] $=((\mu_{A}\circ((\lambda_{A}\ast id_{A})\otimes A))\otimes(\mu_{H}\circ(\lambda_{H}\otimes \mu_{H})))\circ(A\otimes((c_{H,A}\otimes H)\circ(H\otimes\Psi)\circ(\delta_{H}\otimes A))\otimes H)$ {\footnotesize\textnormal{(by}}
\item[] {\footnotesize\textnormal{naturality of $c$, the symmetric character of {\sf C} and associativity of $\mu_{A}$)}}
\item[]$=\varepsilon_{A}\otimes ((A\otimes (\mu_{H}\circ (\lambda_{H}\otimes(\mu_{H}\circ (\phi_{H}\otimes H)))))\circ (((\varphi_{A}\otimes H)\circ (H\otimes c_{H,A})\circ((c_{H,H}\circ\delta_{H})\otimes A))\otimes H\otimes $
\item [] $A\otimes H)\circ (((H\otimes c_{H,A}\otimes A)\circ(\delta_{H}\otimes\delta_{A}))\otimes H))$ {\footnotesize\textnormal{(by \eqref{antipode}, naturality of $c$ and coassociativity of $\delta_{H}$)}}
\item[] $= \varepsilon_{A}\otimes((\varphi_{A}\otimes(\mu_{H}\circ(\overline{\Omega}\otimes H)))\circ(\delta_{H\otimes A}\otimes H))$ {\footnotesize\textnormal{(by \eqref{cocom_class}, coassociativity of $\delta_{H}$, naturality of $c$ and associativity}}
\item[] {\footnotesize\textnormal{of $\mu_{H}$)}}.
\end{itemize}
As a consequence, the following is obtained:
\begin{itemize}
\itemindent=-32pt 
\item[ ]$\hspace{0.38cm}\mu_{A\otimes H}\circ(\mu_{A \bowtie H}\otimes\Gamma_{A\otimes H})\circ(A\otimes H\otimes c_{A\otimes H,A\otimes H}\otimes A\otimes H)\circ (\delta_{A\otimes H}\otimes A\otimes H\otimes A\otimes H)$
\item[] $=(\mu_{A}\otimes\mu_{H})\circ (A\otimes c_{H,A}\otimes H)\circ(((\mu_{A}\otimes\mu_{H})\circ(A\otimes \Psi\otimes H))\otimes((\varphi_{A}\otimes(\mu_{H}\circ (\overline{\Omega}\otimes H)))\circ(\delta_{H\otimes A}\otimes$
\item[] $H)))\circ(A\otimes((H\otimes A\otimes c_{H,H})\circ(H\otimes c_{H,A}\otimes H)\circ(\delta_{H}\otimes A\otimes H))\otimes A\otimes H)$ {\footnotesize\textnormal{(by naturality of $c$ and counit}}
\item[] {\footnotesize\textnormal{property for $A$)}}
\item[] $=((\mu_{A}\circ((\mu_{A}\circ(A\otimes\varphi_{A}))\otimes A))\otimes(\mu_{H}\circ ((\mu_{H}\circ (\phi_{H}\otimes H))\otimes(\mu_{H}\circ (\overline{\Omega}\otimes H)))))\circ(A\otimes H\otimes A\otimes$
\item[] $((\varphi_{A}\otimes H)\circ(H\otimes c_{H,A})\circ((c_{H,H}\circ\delta_{H})\otimes A))\otimes A\otimes H\otimes H\otimes A\otimes H)\circ (A\otimes ((H\otimes c_{H,A})\circ(\delta_{H}\otimes $
\item[] $A))\otimes((c_{A,A}\otimes c_{H,H})\circ (A\otimes c_{H,A}\otimes H)\circ(c_{H,A}\otimes c_{H,A}))\otimes A\otimes H)\circ(A\otimes\delta_{H\otimes A}\otimes H\otimes\delta_{A}\otimes H)$ {\footnotesize\textnormal{(by}}
\item[] {\footnotesize\textnormal{naturality of $c$ and the symmetric character of {\sf C})}}
\item[] $=(\mu_{A}\otimes\mu_{H})\circ (A\otimes(\mu_{A}\circ(\varphi_{A}\otimes\varphi_{A})\circ(H\otimes c_{H,A}\otimes A)\circ(\delta_{H}\otimes A\otimes A))\otimes(\mu_{H}\circ(\mu_{H}\otimes \overline{\Omega})\circ(H\otimes$
\item[] $c_{H,H}\otimes A)\circ(((\phi_{H}\otimes H)\circ(H\otimes c_{H,A})\circ(\delta_{H}\otimes A))\otimes H\otimes A))\otimes H)\circ(A\otimes H\otimes A\otimes((c_{H,A}\otimes A\otimes H)\circ$
\item[] $(H\otimes c_{A,A}\otimes H)\circ(H\otimes A\otimes c_{H,A}))\otimes A\otimes H)\circ(A\otimes \delta_{H\otimes A}\otimes H\otimes\delta_{A}\otimes H)$ {\footnotesize\textnormal{(by \eqref{cocom_class}, coassociativity of}}
\item[] {\footnotesize\textnormal{$\delta_{H}$ and naturality of $c$)}}
\item[] $=(\mu_{A}\otimes\mu_{H})\circ(A\otimes(\varphi_{A}\circ(H\otimes\mu_{A}))\otimes(\mu_{H}\circ((\phi_{H}\circ(H\otimes\mu_{A}))\otimes H)\circ(H\otimes A\otimes c_{H,A}))\otimes H)\circ(A\otimes H\otimes $
\item[] $A\otimes((c_{H,A}\otimes A\otimes H)\circ(H\otimes c_{A,A}\otimes H)\circ(H\otimes A\otimes c_{H,A}))\otimes A\otimes H)\circ(A\otimes\delta_{H\otimes A}\otimes H\otimes\delta_{A}\otimes H)$ 
\item[] {\footnotesize\textnormal{(by \eqref{mod_alg_z1} and \eqref{z2})}}
\item[] $=(\mu_{A}\otimes\mu_{H})\circ(A\otimes((\varphi_{A}\otimes\phi_{H})\circ(H\otimes c_{H,A}\otimes A)\circ(\delta_{H}\otimes ((\mu_{A}\otimes\mu_{A})\circ\delta_{A\otimes A})))\otimes\mu_{H})\circ(A\otimes H\otimes A\otimes c_{H,A}\otimes H)$ 
\item[] {\footnotesize\textnormal{(by naturality of $c$ and associativity of $\mu_{H}$)}}
\item[] $=\mu_{A\bowtie H}\circ(A\otimes H\otimes\mu_{A\otimes H})$ {\footnotesize\textnormal{(by the condition of coalgebra morphism for $\mu_{A}$)}},
\end{itemize}
i.e., $(A\otimes H,A\bowtie H)$ is a Hopf brace in {\sf C}.

On the other hand, to obtain \eqref{mod_alg_z1} from \eqref{compatHbrace} for $(A\otimes H,A\bowtie H)$ it is enough to compose on the right with $\eta_{A}\otimes H\otimes A\otimes\eta_{H}\otimes A\otimes\eta_{H}$, and on the left with $A\otimes\varepsilon_{H}$, whereas to obtain \eqref{z2} from \eqref{compatHbrace} for $(A\otimes H,A\bowtie H)$ we only have to compose on the right with $\eta_{A}\otimes H\otimes A\otimes H\otimes A\otimes\eta_{H}$, and on the left with $\varepsilon_{A}\otimes H$.
\end{proof}
\begin{corollary}\label{cocom_bicros}
Let $A$ and $H$ be Hopf algebras in {\sf C}. If $(A,\varphi_{A})$ is a left $H$-module algebra-coalgebra such that \eqref{cocom_class} holds, then $(A\otimes H,A\sharp H)$ is a Hopf brace in {\sf C}.
\end{corollary}
\begin{proof} 
The fact that $(A,\varphi_{A})$ is a left $H$-module algebra-coalgebra satisfying \eqref{cocom_class} is equivalent to the fact that $(A,H,\varphi_{A},\phi_{H}\coloneqq H\otimes\varepsilon_{A})$ is a matched pair of Hopf algebras. Thus, we can apply the previous theorem to this particular situation resulting the following: On the one hand, condition \eqref{mod_alg_z1} automatically holds because $(A,\varphi_{A})$ is a left $H$-module algebra and, on the other hand, condition \eqref{z2} follows by taking into account that $\overline{\Omega}=(\eta_{H}\circ\varepsilon_{H})\otimes\varepsilon_{A}$ in this case. Then, due to the fact that the bicrossed product Hopf algebra built from the matched pair $(A,H,\varphi_{A},\phi_{H}\coloneqq H\otimes\varepsilon_{A})$ is the smash product Hopf algebra $A\sharp H$, we obtain that $(A\otimes H,A\sharp H)$ is a Hopf brace.
\end{proof}
\begin{remark}
The construction of the previous corollary had already been studied in \cite[Example 1.5]{AGV}, but in this example it is required both Hopf algebras, $A$ and $H$, to be cocommutative. In the above result the hypotheses have been reduced with respect to \cite[Example 1.5]{AGV}: Now there is no requirements on $A$ and it is only demanded that $(A,\varphi_{A})$ is in the cocommutativity class of $H$, which is always true in case that $H$ is cocommutative.
\end{remark} 

In \cite{TD}, Y. Doi and M. Takeuchi proved that the Drinfeld's Double Hopf algebra can be constructed as a bicrossed product Hopf algebra built from a particular matched pair. Thus, Theorem \ref{class_cocom_bicros} and Corollary \ref{cocom_bicros} can be applied in order to obtain conditions under which the tensor product Hopf algebra and the Drinfeld's Double give rise to a Hopf brace. 

Let's recall the construction of the Drinfeld's Double Hopf algebra. Let $H$ be a finite Hopf algebra in {\sf C}, i.e., there exists $H^{\ast}\in{\sf C}$ such that $H\otimes - \dashv H^{\ast}\otimes -$ is a {\sf C}-adjunction with unit $a_{H}$ given by the family of natural morphisms
\[\{a_{H}(X)\colon X\rightarrow H^{\ast}\otimes H\otimes X\}_{X\in{\sf C}},\]
and counit $b_{H}$ given by
\[\{b_{H}(X)\colon H\otimes H^{\ast}\otimes X\rightarrow X\}_{X\in{\sf C}}.\] 

In the finiteness case, it is well-known that $\lambda_{H}$ is always an isomorphism in {\sf C} (see \cite[Corollary 1]{LopVillan}). Then, $\widehat{H}\coloneqq (H^{cop})^{\ast}$ is a Hopf algebra with
\begin{gather*}
\eta_{\widehat{H}}\coloneqq (H^{\ast}\otimes \varepsilon_{H})\circ a_{H}(K),\\ \mu_{\wh}\coloneqq (H^{\ast}\otimes b_{H}(K))\circ(H^{\ast}\otimes H\otimes b_{H}(K)\otimes H^{\ast})\circ(((H^{\ast}\otimes (c_{H,H}\circ\delta_{H}))\circ a_{H}(K))\otimes H^{\ast}\otimes H^{\ast}),\\
\varepsilon_{\wh}\coloneqq b_{H}(K)\circ(\eta_{H}\otimes H^{\ast}),\\
\delta_{\wh}\coloneqq (H^{\ast}\otimes H^{\ast}\otimes(b_{H}(K)\circ (\mu_{H}\otimes H^{\ast})))\circ(((H^{\ast}\otimes a_{H}(K)\otimes H)\circ a_{H}(K))\otimes H^{\ast}),\\
\lambda_{\wh}\coloneqq (H^{\ast}\otimes b_{H}(K))\circ(H^{\ast}\otimes\lambda_{H}^{-1}\otimes H^{\ast})\circ(a_{H}(K)\otimes H^{\ast}).
\end{gather*}

Let's denote by $T(H)\coloneqq \wh\otimes H$, the usual tensor product Hopf algebra of $\wh$ with $H$. Moreover, it is well known (see \cite{AFG_2coc_HCuasi} for details) that, given a Hopf algebra $A=(A,\eta_{A},\mu_{A},\varepsilon_{A},\delta_{A},\lambda_{A})$ in {\sf C}, if $\sigma\colon A\otimes A\rightarrow K$ is an invertible normal  2-cocycle on $A$, then it is possible to construct a new Hopf algebra structure on $A$ which shares the same coalgebra structure and the same unit with $A$, but whose product and antipode are defined as follows:
\begin{gather}
\label{prod_deform_2_coc}\mu_{A}^{\sigma}\coloneqq(\sigma\otimes\mu_{A}\otimes\sigma^{-1})\circ(A\otimes A\otimes\delta_{A\otimes A})\circ\delta_{A\otimes A},\\
\label{ant_deform_2_coc}\lambda_{A}^{\sigma}\coloneqq (f\otimes((\lambda_{A}\otimes f^{-1})\circ\delta_{A}))\circ\delta_{A},
\end{gather}
where $\sigma^{-1}$ is the inverse of $\sigma$ for the convolution in $\Hom_{{\sf C}}(A\otimes A,K)$, $f\coloneqq \sigma\circ(A\otimes\lambda_{A})\circ\delta_{A}$ and $f^{-1}$ is the inverse of $f$ for the convolution in $\Hom_{{\sf C}}(A,K)$, which coincides in this situation with $f^{-1}=\sigma^{-1}\circ(\lambda_{A}\otimes A)\circ\delta_{A}$. This Hopf algebra structure is the deformation of $A$ by $\sigma$ and it will be denoted by $A^{\sigma}$. So, in \cite[Proposition 1.5]{TD}, \cite[Remark 2.3]{TD} and \cite[Remark 2.6]{TD}, Y. Doi and M. Takeuchi assert that the morphism 
\begin{equation}\label{2coc_D(H)}
\omega\coloneqq \varepsilon_{\wh}\otimes b_{H}(K)\otimes \varepsilon_{H}
\end{equation}
is an invertible normal 2-cocycle on $T(H)$ with inverse \begin{equation}\label{2coc_D(H)_inverse}
	\omega^{-1}= \varepsilon_{\wh}\otimes(b_{H}(K)\circ(\lambda_{H}^{-1}\otimes H^{\ast}))\otimes\varepsilon_{H},
\end{equation}  
and the Drinfeld's Double Hopf algebra, $D(H)$, is equal to the deformation of $T(H)$ by $\omega$, i.e., 
\[D(H)=T(H)^{\omega}.\]

Following formulas \eqref{prod_deform_2_coc} and \eqref{ant_deform_2_coc} for the product and the antipode, respectively, together with the expressions \eqref{2coc_D(H)} and \eqref{2coc_D(H)_inverse} for the 2-cocycle $\omega$ and its inverse, it is easy to show that 
\begin{equation}\label{for_prod_D(H)}\mu_{D(H)}=\mu_{T(H)}\circ(H^{\ast}\otimes((b_{H}(K)\otimes H\otimes H^{\ast})\circ(\delta_{H\otimes\wh}\otimes(b_{H}(K)\circ(\lambda_{H}^{-1}\otimes H^{\ast})))\circ\delta_{H\otimes\wh})\otimes H),\end{equation}
and 
\begin{equation*}\label{for_ant_D(H)}
\lambda_{D(H)}=((b_{H}(K)\circ (H\otimes\lambda_{\wh})\circ c_{H^{\ast},H})\otimes((\lambda_{\wh}\otimes\lambda_{H}\otimes(b_{H}(K)\circ c_{H^{\ast},H}))\circ \delta_{\wh\otimes H}))\circ\delta_{\wh\otimes H}.
\end{equation*}

Furthermore, Y. Doi and M. Takeuchi also prove that $D(H)$ is equal to the bicrossed product Hopf algebra $\wh\bowtie H$ built from the matched pair $(\wh,H,\varphi_{\wh},$ $\phi_{H})$, where
\begin{gather}
\label{varphi_wh}\varphi_{\wh}\coloneqq (H^{\ast}\otimes b_{H}(K))\circ(H^{\ast}\otimes R\otimes H^{\ast})\circ(a_{H}(K)\otimes H\otimes H^{\ast}),\\
\label{phi_H}\phi_{H}\coloneqq (H\otimes b_{H}(K))\circ(\overline{R}\otimes H^{\ast}),
\end{gather}
and $R$ and $\overline{R}$ are defined as follows:
\begin{gather}
\label{R_def}R\coloneqq \mu_{H}\circ c_{H,H}\circ (\mu_{H}\otimes\lambda_{H}^{-1})\circ(H\otimes\delta_{H}),\\
\label{Rbarra_def}\overline{R}\coloneqq (H\otimes(\mu_{H}\circ(\lambda_{H}^{-1}\otimes H)))\circ(\delta_{H}\otimes H)\circ c_{H,H}\circ\delta_{H}.
\end{gather}

As a consequence, 
\[D(H)=T(H)^{\omega}=\wh\bowtie H.\]

Then, the following results about under what conditions $(T(H),D(H))$ is a Hopf brace in {\sf C} are obtained:
\begin{theorem}\label{D(H)_cocommutative_HBr}
Let $H$ be a finite and cocommutative Hopf algebra in {\sf C}. Then, $(T(H),D(H))$ is a Hopf brace in {\sf C}.
\end{theorem}
\begin{proof}
If $H$ is cocommutative, then it is straightforward to show that $\phi_{H}=H\otimes \varepsilon_{\wh}$, and then $D(H)=\wh\sharp H$. Thus, $(T(H), D(H))$ is a Hopf brace in {\sf C} by Corollary \ref{cocom_bicros}.
\end{proof}

\begin{example}\label{D(F[G])_cocom}
Assume that ${\sf C}={}_{\mathbb{F}}{\sf Vect}$, the category of vector spaces over the field $\mathbb{F}$. Let $(G,\cdot)$ be a finite group with unit $e$. Then, the group algebra $\mathbb{F}[G]=\bigoplus_{g\in G}\mathbb{F}g$ is a cocommutative Hopf algebra in ${}_{\mathbb{F}}{\sf Vect}$ with
\begin{gather*}
\eta_{\mathbb{F}[G]}(1_{\mathbb{F}})=e,\quad \mu_{\mathbb{F}[G]}(g\otimes h)=g\cdot h,\quad
\varepsilon_{\mathbb{F}[G]}(g)=1_{\mathbb{F}},\quad \delta_{\mathbb{F}[G]}(g)=g\otimes g,\quad
\lambda_{\mathbb{F}[G]}(g)=g^{-1}
\end{gather*}
for all $g,h\in G$. 

Note that, due to the cocommutativity of $\mathbb{F}[G]$, $\widehat{\mathbb{F}[G]}=\mathbb{F}[G]^{\ast}$, which is a commutative Hopf algebra in ${}_{\mathbb{F}}{\sf Vect}$ with
\begin{gather*}
\eta_{\mathbb{F}[G]^{\ast}}(1_{\mathbb{F}})=\sum_{r\in G}f_{r},\quad \mu_{\mathbb{F}[G]^{\ast}}(f_{r}\otimes f_{s})=\delta_{r,s}f_{r},\\
\varepsilon_{\mathbb{F}[G]^{\ast}}(f_{r})=\delta_{r,e},\quad \delta_{\mathbb{F}[G]^{\ast}}(f_{r})=\sum_{\substack{l,m\in G\\r=l\cdot m}}f_{l}\otimes f_{m},\quad
\lambda_{\mathbb{F}[G]^{\ast}}(f_{r})=f_{r^{-1}},
\end{gather*}
where, if $\{g\}_{g\in G}$ is a basis of $\mathbb{F}[G]$, $\{f_{r}\}_{r\in G}$ is its respective dual basis, i.e., $f_{r}(g)=\delta_{r,g}$.

By the previous theorem, $(T(\mathbb{F}[G]),D(\mathbb{F}[G]))$ is a Hopf brace in {\sf C}, being $T(\mathbb{F}[G])=\mathbb{F}[G]^{\ast}\otimes \mathbb{F}[G]$, which is a Hopf algebra with
\begin{gather*}
\eta_{T(\mathbb{F}[G])}(1_{\mathbb{F}})=\sum_{r\in G}f_{r}\otimes e,\quad \mu_{T(\mathbb{F}[G])}(f_{r}\otimes g\otimes f_{s}\otimes h)=\delta_{r,s}f_{r}\otimes g\cdot h,\\
\varepsilon_{T(\mathbb{F}[G])}(f_{r}\otimes g)=\delta_{r,e},\quad \delta_{T(\mathbb{F}[G])}(f_{r}\otimes g)=\sum_{\substack{l,m\in G\\r=l\cdot m}}f_{l}\otimes g\otimes f_{m}\otimes g,\quad
\lambda_{T(\mathbb{F}[G])}(f_{r}\otimes g)=f_{r^{-1}}\otimes g^{-1},
\end{gather*}
whereas $D(\mathbb{F}[G])=\mathbb{F}[G]^{\ast}\sharp \mathbb{F}[G]$, which is a Hopf algebra with
\begin{gather*}
\eta_{D(\mathbb{F}[G])}=\eta_{T(\mathbb{F}[G])},\quad \varepsilon_{D(\mathbb{F}[G])}=\varepsilon_{T(\mathbb{F}[G])},\quad \delta_{D(\mathbb{F}[G])}=\delta_{T(\mathbb{F}[G])},\\
\mu_{D(\mathbb{F}[G])}(f_{r}\otimes g\otimes f_{s}\otimes h)=\delta_{r,g^{-1}\cdot s\cdot g}f_{r}\otimes g\cdot h,\quad\lambda_{D(\mathbb{F}[G])}(f_{r}\otimes g)=f_{g\cdot r^{-1}\cdot g^{-1}}\otimes g^{-1}.
\end{gather*}
In this case, $\mathbb{F}[G]^{\ast}$ is a left $\mathbb{F}[G]$-module algebra-coalgebra with 
\begin{gather*}
\varphi_{\mathbb{F}[G]^{\ast}}(g\otimes f_{r})=f_{g\cdot r\cdot g^{-1}} .
\end{gather*}

Following \cite[Example 7.1.9]{MAJ_Book}, $D(\mathbb{F}[G])$ can be interpreted as the quantum algebra of observables of a quantum system.

If we also require $(G,\cdot)$ to be commutative, then $\mathbb{F}[G]$ is a commutative Hopf algebra and, as a consequence, $\mathbb{F}[G]^{\ast}$ is a cocommutative Hopf algebra in ${}_{\mathbb{F}}{\sf Vect}$. So, $(T(\mathbb{F}[G]^{\ast}),D(\mathbb{F}[G]^{\ast}))$ is a Hopf brace in ${}_{\mathbb{F}}{\sf Vect}$ too by applying Theorem \ref{D(H)_cocommutative_HBr}.
\end{example}

In addition, we can also subtract from Theorem \ref{class_cocom_bicros} what happens when $(\wh,\varphi_{\wh})$ belongs to the cocommutativity class of $H$ rather than $H$ being cocommutative, which is a weaker condition.
\begin{corollary}\label{D(H)_cocomclass_HBr}
Let $H$ be a finite Hopf algebra in {\sf C}. If $H$ satisfies that 
\begin{equation}\label{cocom_class_D(H)}
(R\otimes H)\circ(H\otimes(c_{H,H}\circ\delta_{H}))=(R\otimes H)\circ(H\otimes\delta_{H}),
\end{equation}
where $R$ is the morphism defined in \eqref{R_def}, then $(T(H),D(H))$ is a Hopf brace in {\sf C} if and only if $R$ is a coalgebra morphism and the equality
\begin{align}\label{z2_D(H)}
&(\mu_{H}\otimes H\otimes H)\circ(H\otimes c_{H,H}\otimes H)\circ(\overline{R}\otimes((\mu_{H}\otimes H)\circ(H\otimes S)\circ c_{H,H}))\circ(\delta_{H}\otimes H)\\\nonumber=&(\mu_{H}\circ\delta_{H})\circ(H\otimes c_{H,H})\circ(\overline{R}\otimes H)
\end{align}
holds, where $S$ is given by 
\begin{equation}\label{S_def}
S\coloneqq (\mu_{H}\otimes H)\circ(\lambda_{H}\otimes\overline{R})\circ\delta_{H},
\end{equation}
and $\overline{R}$ is the morphism defined in \eqref{Rbarra_def}.
\end{corollary}
\begin{proof}
This is a particular case of Theorem \ref{class_cocom_bicros} by considering the matched pair $(\wh,H,\varphi_{\wh},\phi_{H})$. Note first that, when $H$ is a finite object in the symmetric monoidal category {\sf C}, the equalities
\begin{gather}\label{fin_unit}
(c_{H,H^{\ast}}\otimes H)\circ(H\otimes a_{H}(K))=(H^{\ast}\otimes c_{H,H})\circ(a_{H}(K)\otimes H),\\
\label{fin_counit}
(b_{H}(K)\otimes H)\circ(H\otimes c_{H,H^{\ast}})=(H\otimes b_{H}(K))\circ(c_{H,H}\otimes H^{\ast})
\end{gather}
hold, which will be very useful in the following computations (see \cite[Lemma 1.16]{FGR} for details).

Under these conditions, \eqref{cocom_class_D(H)} is equivalent to \eqref{cocom_class} for $\varphi_{\wh}$. Indeed, by \eqref{varphi_wh}, \eqref{cocom_class} can be written equivalently as follows:
\begin{align}\label{p1}
&(H^{\ast}\otimes((b_{H}(K)\otimes H)\circ (R\otimes c_{H,H^{\ast}})))\circ(a_{H}(K)\otimes(c_{H,H}\circ\delta_{H})\otimes H^{\ast})\\\nonumber=&(H^{\ast}\otimes((b_{H}(K)\otimes H)\circ (R\otimes c_{H,H^{\ast}})))\circ(a_{H}(K)\otimes\delta_{H}\otimes H^{\ast}).
\end{align}

Then, by using \eqref{fin_counit}, \eqref{p1} is equivalent to 
\begin{align}\label{p2}
&(H^{\ast}\otimes H\otimes b_{H}(K))\circ(H^{\ast}\otimes(c_{H,H}\circ(R\otimes H))\otimes H^{\ast})\circ(a_{H}(K)\otimes(c_{H,H}\circ\delta_{H})\otimes H^{\ast})\\\nonumber=&(H^{\ast}\otimes H\otimes b_{H}(K))\circ(H^{\ast}\otimes(c_{H,H}\circ(R\otimes H))\otimes H^{\ast})\circ(a_{H}(K)\otimes\delta_{H}\otimes H^{\ast}).
\end{align}

Moreover, by the properties of the unit and the counit of the {\sf C}-adjunction $H\otimes - \dashv H^{\ast}\otimes -$, \eqref{p2} is equivalent to 
\begin{align*}
c_{H,H}\circ(R\otimes H)\circ(H\otimes (c_{H,H}\circ\delta_{H}))=c_{H,H}\circ(R\otimes H)\circ(H\otimes \delta_{H}),
\end{align*}
which is equivalent to \eqref{cocom_class_D(H)} because $c_{H,H}$ is an isomorphism.

So, by Theorem \ref{class_cocom_bicros}, if \eqref{cocom_class_D(H)} holds, then $(T(H),D(H))$ is a Hopf brace if and only if \eqref{mod_alg_z1} and \eqref{z2} hold. Let's see that \eqref{mod_alg_z1} holds if and only if $R$ is a coalgebra morphism. On the one hand, the right hand side of equality \eqref{mod_alg_z1} becomes
\begin{itemize}
\itemindent=-32pt
\item[ ]$\hspace{0.38cm}\varphi_{\wh}\circ (H\otimes\mu_{\wh})$
\item[] $=(H^{\ast}\otimes(b_{H}(K)\circ(H\otimes b_{H}(K)\otimes H^{\ast})))\circ(((H^{\ast}\otimes(c_{H,H}\circ\delta_{H}\circ R))\circ(a_{H}(K)\otimes H))\otimes H^{\ast}\otimes H^{\ast})$
\end{itemize}
by using \eqref{varphi_wh} and the properties of the unit and the counit of the {\sf C}-adjunction $H\otimes - \dashv H^{\ast}\otimes -$, whereas, on the other hand, the left hand side is expanded as follows:
\begin{itemize}
\itemindent=-32pt 
\item[ ]$\hspace{0.38cm}\mu_{\wh}\circ(\varphi_{\wh}\otimes\varphi_{\wh})\circ(H\otimes c_{H,H^{\ast}}\otimes H^{\ast})\circ (\delta_{H}\otimes H^{\ast}\otimes H^{\ast})$
\item[] $=(H^{\ast}\otimes(b_{H}(K)\circ (R\otimes H^{\ast})))\circ(H^{\ast}\otimes H\otimes (b_{H}(K)\circ (R\otimes H^{\ast}))\otimes H\otimes H^{\ast})\circ(((H^{\ast}\otimes(c_{H,H}\circ\delta_{H}))\circ $
\item[] $a_{H}(K))\otimes((H\otimes c_{H,H^{\ast}})\circ(\delta_{H}\otimes H^{\ast}))\otimes H^{\ast})$ {\footnotesize\textnormal{(by \eqref{varphi_wh} and the properties of the unit and the counit of the}}
\item[] {\footnotesize\textnormal{adjunction)}}
\item[] $=(H^{\ast}\otimes(b_{H}(K)\circ (H\otimes b_{H}(K)\otimes H^{\ast})))\circ(((H^{\ast}\otimes((R\otimes R)\circ (H\otimes c_{H,H}\otimes H)\circ ((c_{H,H}\circ\delta_{H})\otimes(c_{H,H}\circ$
\item[] $\delta_{H}))))\circ(a_{H}(K)\otimes H))\otimes H^{\ast}\otimes H^{\ast})$ {\footnotesize\textnormal{(by \eqref{fin_counit} and naturality of $c$)}}.
\end{itemize}
Then, equality \eqref{mod_alg_z1} holds if and only if 
\begin{itemize}
\itemindent=-32pt 
\item[ ]$\hspace{0.38cm}(H^{\ast}\otimes(b_{H}(K)\circ(H\otimes b_{H}(K)\otimes H^{\ast})))\circ(((H^{\ast}\otimes(c_{H,H}\circ\delta_{H}\circ R))\circ(a_{H}(K)\otimes H))\otimes H^{\ast}\otimes H^{\ast})$
\item[] $=(H^{\ast}\otimes(b_{H}(K)\circ (H\otimes b_{H}(K)\otimes H^{\ast})))\circ(((H^{\ast}\otimes((R\otimes R)\circ (H\otimes c_{H,H}\otimes H)\circ ((c_{H,H}\circ\delta_{H})\otimes(c_{H,H}\circ$
\item[] $\delta_{H}))))\circ(a_{H}(K)\otimes H))\otimes H^{\ast}\otimes H^{\ast})$,
\end{itemize}
which is equivalent to
\begin{align}\label{p3}
c_{H,H}\circ\delta_{H}\circ R=(R\otimes R)\circ(H\otimes c_{H,H}\otimes H)\circ((c_{H,H}\circ\delta_{H})\otimes(c_{H,H}\circ\delta_{H}))
\end{align}
by the properties of the unit and the counit of the {\sf C}-adjunction $H\otimes - \dashv H^{\ast}\otimes -$. So, composing on the left with the isomorphism $c_{H,H}$, \eqref{p3} is equivalent to the fact that $R$ is a coalgebra morphism.

To conclude it only remains to prove that \eqref{z2} is equivalent to \eqref{z2_D(H)}. In this case, note that the morphism $\overline{\Omega}$ defined in \eqref{Omegabar_def} satisfies that
\begin{align}\label{p4}
&\overline{\Omega}
\\\nonumber=&\mu_{H}\circ(\lambda_{H}\otimes((H\otimes b_{H}(K))\circ (\overline{R}\otimes H^{\ast})))\circ(\delta_{H}\otimes H^{\ast})\;\footnotesize\textnormal{(by \eqref{phi_H})}\\\nonumber=&(H\otimes b_{H}(K))\circ(S\otimes H^{\ast})\;\footnotesize\textnormal{(by \eqref{S_def})}.
\end{align}

Then, the right hand side of \eqref{z2} becomes 
\begin{itemize}
\itemindent=-32pt 
\item[ ]$\hspace{0.38cm}\mu_{H}\circ((\phi_{H}\circ(H\otimes\mu_{\wh}))\otimes H)\circ(H\otimes H^{\ast}\otimes c_{H,H^{\ast}})$ 
\item[] $=\mu_{H}\circ(H\otimes (b_{H}(K)\circ (H\otimes b_{H}(K)\otimes H^{\ast}))\otimes H)\circ(((H\otimes(c_{H,H}\circ\delta_{H}))\circ\overline{R})\otimes H^{\ast}\otimes c_{H,H^{\ast}})$ {\footnotesize\textnormal{(by \eqref{phi_H}}}
\item[] {\footnotesize\textnormal{and the properties of the unit and the counit of the adjunction)}}
\item[] $=(\mu_{H}\otimes b_{H}(K))\circ (H\otimes c_{H,H}\otimes H^{\ast})\circ (((H\otimes H\otimes b_{H}(K))\circ (((H\otimes(c_{H,H}\circ\delta_{H}))\circ\overline{R})\otimes H^{\ast}))\otimes H\otimes H^{\ast})$ 
\item[] {\footnotesize\textnormal{(by \eqref{fin_counit})}},
\end{itemize}
whereas the left hand side of \eqref{z2} is given by 
\begin{itemize}
\itemindent=-32pt 
\item[ ]$\hspace{0.38cm}\mu_{H}\circ(\mu_{H}\otimes\overline{\Omega})\circ(H\otimes c_{H,H}\otimes H^{\ast})\circ(((\phi_{H}\otimes H)\circ(H\otimes c_{H,H^{\ast}})\circ(\delta_{H}\otimes H^{\ast}))\otimes H\otimes H^{\ast})$
\item[] $=\mu_{H}\circ ((\mu_{H}\circ (H\otimes b_{H}(K)\otimes H)\circ (\overline{R}\otimes H^{\ast}\otimes H))\otimes((H\otimes b_{H}(K))\circ (S\otimes H^{\ast})))\circ (((H\otimes H^{\ast}\otimes$ 
\item[] $c_{H,H})\circ(H\otimes c_{H,H^{\ast}}\otimes H)\circ(\delta_{H}\otimes H^{\ast}\otimes H))\otimes H^{\ast})$ {\footnotesize\textnormal{(by \eqref{phi_H}, \eqref{p4} and the properties of the unit and the}}
\item[] {\footnotesize\textnormal{counit of the adjunction)}}
\item[] $=\mu_{H}\circ(\mu_{H}\otimes((H\otimes b_{H}(K))\circ(S\otimes H^{\ast})))\circ(H\otimes c_{H,H}\otimes H^{\ast})\circ(((H\otimes H\otimes b_{H}(K))\circ(((H\otimes c_{H,H})\circ$
\item[] $(\overline{R}\otimes H)\circ\delta_{H})\otimes H^{\ast}))\otimes H\otimes H^{\ast})$ {\footnotesize\textnormal{(by \eqref{fin_counit})}}.
\end{itemize}

As a consequence, \eqref{z2} is verified if and only if the following equality holds:
\begin{itemize}
\itemindent=-32pt 
\item[ ]$\hspace{0.38cm}(\mu_{H}\otimes b_{H}(K))\circ (H\otimes c_{H,H}\otimes H^{\ast})\circ (((H\otimes H\otimes b_{H}(K))\circ (((H\otimes(c_{H,H}\circ\delta_{H}))\circ\overline{R})\otimes H^{\ast}))\otimes H\otimes H^{\ast})$
\item[] $=\mu_{H}\circ(\mu_{H}\otimes((H\otimes b_{H}(K))\circ(S\otimes H^{\ast})))\circ(H\otimes c_{H,H}\otimes H^{\ast})\circ(((H\otimes H\otimes b_{H}(K))\circ(((H\otimes c_{H,H})\circ$
\item[] $(\overline{R}\otimes H)\circ\delta_{H})\otimes H^{\ast}))\otimes H\otimes H^{\ast})$,
\end{itemize}
which is equivalent to
\begin{align}\label{p5}
&(\mu_{H}\otimes H)\circ(H\otimes c_{H,H})\circ(H\otimes((H\otimes b_{H}(K))\circ((c_{H,H}\circ\delta_{H})\otimes H^{\ast}))\otimes H)\circ(\overline{R}\otimes c_{H,H^{\ast}})\\\nonumber=& (\mu_{H}\otimes H)\circ(\mu_{H}\otimes S)\circ(H\otimes c_{H,H})\circ(((H\otimes H\otimes b_{H}(K))\circ(H\otimes c_{H,H}\otimes H^{\ast})\circ(\overline{R}\otimes H \\\nonumber &\hspace{-.38 cm}\otimes H^{\ast}))\otimes H)\circ(\delta_{H}\otimes c_{H,H^{\ast}})
\end{align}
by the properties of the unit and the counit of the {\sf C}-adjunction $H\otimes - \dashv H^{\ast}\otimes -$ and by the fact that $c_{H,H^{\ast}}$ is an isomorphism in {\sf C}. Thus, by applying \eqref{fin_counit}, the properties of the unit and the counit of the {\sf C}-adjunction $H\otimes - \dashv H^{\ast}\otimes -$ and the naturality of $c$, \eqref{p5} is equivalent to
\begin{align*}
&(\mu_{H}\otimes(c_{H,H}\circ\delta_{H}))\circ(H\otimes c_{H,H})\circ(\delta_{H}\otimes H)\\
=&(\mu_{H}\otimes c_{H,H})\circ(H\otimes c_{H,H}\otimes H)\circ(\overline{R}\otimes((\mu_{H}\otimes H)\circ(H\otimes S)\circ c_{H,H}))\circ(\delta_{H}\otimes H)
\end{align*} 
which is equivalent to \eqref{z2_D(H)} by the fact that $c_{H,H}$ is an isomorphism in {\sf C}.
\end{proof}

So far in Theorem \ref{mainth2} we have obtained conditions under which, given Hopf braces $\mathbb{A}=(A_{1},A_{2})$ and $\mathbb{H}=(H_{1},H_{2})$, $(A_{1}\otimes H_{1},A_{2}\bowtie H_{2})$ is a Hopf brace with the bicrossed product appearing in the second component. The following result examines what the conditions must be satisfied in case that the bicrossed product appears in the first component.

\begin{theorem}\label{mainth} 
Let $\mathbb{A}=(A_{1},A_{2})$ and $\mathbb{H}=(H_{1},H_{2})$ be Hopf braces in {\sf C} such that
\begin{itemize}
\item[(i)] $(A_{1},H_{1},\varphi_{A}^{1},\phi_{H})$ is a matched pair of Hopf algebras,
\item[(ii)] $(A_{2},\varphi_{A}^{2})$ is a left $H_{2}$-module algebra-coalgebra satisfying \eqref{adcond_mp_phiHtriv}.
\end{itemize}

The pair $$(A_{1}\bowtie H_{1},A_{2}\sharp H_{2})$$ is a Hopf brace in {\sf C} if and only if the following conditions hold:
\begin{itemize}
\item[(C1)] $\hspace{.38cm}(A\otimes\mu_{H}^{1})\circ(\Psi^{1}\otimes\mu_{H}^{2})\circ(H\otimes\Psi^{2}\otimes H)\circ(((\Gamma'_{H_{1}}\otimes H)\circ(H\otimes c_{H,H})\circ(\delta_{H}\otimes H))\otimes A\otimes H)\\=(A\otimes \mu_{H}^{2})\circ(\Psi^{2}\otimes\mu_{H}^{1})\circ(H\otimes\Psi^{1}\otimes H),$
\item[(C2)] $(\Gamma_{A_{1}}\otimes H)\circ(A\otimes\Psi^{1})=\Psi^{1}\circ(H\otimes\Gamma_{A_{1}})\circ(c_{A,H}\otimes A),$
\item[(C3)] $\Psi^{2}\circ (H\otimes \mu_{A}^{1})=(\mu_{A}^{1}\otimes H)\circ(A\otimes \Psi^{2})\circ(\Psi^{2}\otimes A),$
\end{itemize}
where 
\begin{gather*}
\Psi^{1}=(\varphi_{A}^{1}\otimes\phi_{H})\circ(H\otimes c_{H,A}\otimes A)\circ(\delta_{H}\otimes\delta_{A}),\\
\Psi^{2}=(\varphi_{A}^{2}\otimes H)\circ(H\otimes c_{H,A})\circ(\delta_{H}\otimes A).
\end{gather*}
\end{theorem}
\begin{proof}
At first we start by proving some relevant equalities that will be used in the development of the proof. Concretely, these equalities are the following:
\begin{gather}\label{equal1}
\Psi^{1}\circ(\eta_{H}\otimes A)=A\otimes\eta_{H},\\
\label{equal2}
\Psi^{2}\circ(\eta_{H}\otimes A)=A\otimes\eta_{H},\\
\label{equal3}
\Psi^{2}\circ(H\otimes\eta_{A})=\eta_{A}\otimes H,\\
\label{equaleps}
(A\otimes\varepsilon_{H})\circ\Psi^{2}=\varphi_{A}^{2},\\
\label{equal4}
(\mu_{A}^{1}\otimes H)\circ (A\otimes \Psi^{1})\circ(\Psi^{1}\otimes A)=\Psi^{1}\circ (H\otimes\mu_{A}^{1}),\\
\label{equal5}
(A\otimes\mu_{H}^{1})\circ (\Psi^{1}\otimes H)\circ(H\otimes\Psi^{1})=\Psi^{1}\circ (\mu_{H}^{1}\otimes A),\\
\label{equal6}
(\Psi^{2}\otimes H)\circ(H\otimes c_{H,A})\circ(\delta_{H}\otimes A)=(A\otimes\delta_{H})\circ\Psi^{2},\\
\label{equal7}
\Gamma_{A_{1}\bowtie H_{1}}=(A\otimes\mu_{H}^{1})\circ((\Psi^{1}\circ (\lambda_{H}^{1}\otimes\Gamma_{A_{1}}))\otimes\mu_{H}^{2})\circ(((c_{A,H}\otimes\Psi^{2})\circ (A\otimes\delta_{H}\otimes A))\otimes H),\\
\label{equal8}
\Gamma_{A_{1}\bowtie H_{1}}\circ(\eta_{A}\otimes H\otimes A\otimes H)=(A\otimes\mu_{H}^{1})\circ (\Psi^{1}\otimes \mu_{H}^{2})\circ (((\lambda_{H}^{1}\otimes\Psi^{2})\circ (\delta_{H}\otimes A))\otimes H),\\
\label{equal9}
\Gamma_{A_{1}\bowtie H_{1}}\circ(\eta_{A}\otimes H\otimes A\otimes \eta_{H})=(A\otimes\mu_{H}^{1})\circ (\Psi^{1}\otimes H)\circ (\lambda_{H}^{1}\otimes\Psi^{2})\circ (\delta_{H}\otimes A),\\
\label{equal10}
\Gamma_{A_{1}\bowtie H_{1}}\circ(A\otimes \eta_{H}\otimes A\otimes \eta_{H})=\Gamma_{A_{1}}\otimes\eta_{H}.
\end{gather}

The proof of \eqref{equal1} follows by the condition of coalgebra morphism for $\eta_{H}$, the naturality of $c$, the axioms of left module for $(A,\varphi_{A}^{1})$, (iii) of Definition \ref{MPdef} and counit property. Using similar arguments, \eqref{equal2} and \eqref{equal3} also hold. Besides, by naturality of $c$ and counit property, \eqref{equaleps} holds. Let's detail the proof of \eqref{equal4}:
\begin{itemize}
\itemindent=-32pt 
\item[ ]$\hspace{0.38cm}(\mu_{A}^{1}\otimes H)\circ (A\otimes \Psi^{1})\circ(\Psi^{1}\otimes A)$
\item[] $=(\mu_{A}^{1}\otimes H)\circ(A\otimes((\varphi_{A}^{1}\otimes\phi_{H})\circ (H\otimes c_{H,A}\otimes A)\circ(((\phi_{H}\otimes\phi_{H})\circ(H\otimes c_{H,A}\otimes A)\circ(\delta_{H}\otimes\delta_{A}))\otimes$
\item[] $\delta_{A})))\circ (((\varphi_{A}^{1}\otimes H\otimes A)\circ(H\otimes c_{H,A}\otimes A)\circ(\delta_{H}\otimes\delta_{A}))\otimes A)$ {\footnotesize\textnormal{(by the condition of coalgebra morphism for}}
\item[] {\footnotesize\textnormal{$\phi_{H}$)}}
\item[] $=((\mu_{A}^{1}\circ (A\otimes\varphi_{A}^{1}))\otimes(\phi_{H}\circ(\phi_{H}\otimes A)))\circ(A\otimes H\otimes ((c_{H,A}\otimes A\otimes A)\circ (H\otimes c_{A,A}\otimes A)\circ(H\otimes A\otimes\delta_{A})))\circ$
\item[] $(((\varphi_{A}^{1}\otimes\phi_{H}\otimes H)\circ(H\otimes c_{H,A}\otimes c_{H,A})\circ(H\otimes H\otimes c_{H,A}\otimes A)\circ(\delta_{H}\otimes H\otimes\delta_{A}))\otimes A\otimes A)\circ(\delta_{H}\otimes\delta_{A}\otimes A)$ 
\item[] {\footnotesize\textnormal{(by naturality of $c$ and coassociativity of $\delta_{A}$ and $\delta_{H}$)}}
\item[] $=((\mu_{A}^{1}\circ (A\otimes \varphi_{A}^{1})\circ (\Psi^{1}\otimes A))\otimes (\phi_{H}\circ (H\otimes\mu_{A}^{1})))\circ (H\otimes((A\otimes c_{H,A}\otimes A)\circ (c_{H,A}\otimes c_{A,A}))\otimes A)\circ(\delta_{H}\otimes\delta_{A}\otimes\delta_{A})$ 
\item[] {\footnotesize\textnormal{(by naturality of $c$ and right module axioms for $(H,\phi_{H})$)}}
\item[] $=((\varphi_{A}^{1}\circ(H\otimes\mu_{A}^{1}))\otimes(\phi_{H}\circ (H\otimes\mu_{A}^{1})))\circ(H\otimes((A\otimes c_{H,A}\otimes A)\circ (c_{H,A}\otimes c_{A,A}))\otimes A)\circ(\delta_{H}\otimes\delta_{A}\otimes\delta_{A})$ 
\item[] {\footnotesize\textnormal{(by (iv) of Definition \ref{MPdef})}}
\item[] $=(\varphi_{A}^{1}\otimes\phi_{H})\circ(H\otimes c_{H,A}\otimes A)\circ(\delta_{H}\otimes((\mu_{A}^{1}\otimes\mu_{A}^{1})\circ(A\otimes c_{A,A}\otimes A)\circ(\delta_{A}\otimes\delta_{A})))$ {\footnotesize\textnormal{(by naturality of $c$)}}
\item[] $=\Psi^{1}\circ(H\otimes\mu_{A}^{1})$ {\footnotesize\textnormal{(by the condition of coalgebra morphism for $\mu_{A}^{1}$)}.}
\end{itemize}

Similarly to the previous one, \eqref{equal5} is also satisfied. The equality \eqref{equal6} is a direct consequence of the coassociativity of $\delta_{H}$ and naturality of $c$. The proof of \eqref{equal7} follows by applying the naturality of $c$ and \eqref{equal4}. Using \eqref{equal7}, the naturality of $c$ and left module axioms for $(A_{1},\Gamma_{A_{1}})$, \eqref{equal8} is obtained. Therefore, by unit property and \eqref{equal8}, \eqref{equal9} holds. To conclude, let's compute \eqref{equal10}:
\begin{itemize}
\itemindent=-32pt 
\item[ ]$\hspace{0.38cm}\Gamma_{A_{1}\bowtie H_{1}}\circ(A\otimes \eta_{H}\otimes A\otimes \eta_{H})$
\item[] $=(A\otimes\mu_{H}^{1})\circ((\Psi^{1}\circ(\lambda_{H}^{1}\otimes\Gamma_{A_{1}}))\otimes \mu_{H}^{2})\circ(\eta_{H}\otimes A\otimes(\Psi^{2}\circ (\eta_{H}\otimes A))\otimes\eta_{H})$ {\footnotesize\textnormal{(by \eqref{equal7}, the condition of }}
\item[] {\footnotesize\textnormal{coalgebra morphism for $\eta_{H}$ and naturality of $c$)}}
\item[] $=\Psi^{1}\circ(\eta_{H}\otimes\Gamma_{A_{1}})$ {\footnotesize\textnormal{(by unit property, \eqref{equal2} and \eqref{u-antip1})}}
\item[] $=\Gamma_{A_{1}}\otimes\eta_{H}$ {\footnotesize\textnormal{(by \eqref{equal1})}.}
\end{itemize}

Note also that,  if (i) is assumed, then it is ensured that $A_{1}\bowtie H_{1}$ is a Hopf algebra in {\sf C}, and also, by Remark \ref{particular_smash}, $A_{2}\sharp H_{2}$ is a Hopf algebra in {\sf C} thanks to (ii).

Taking into account the previous reasoning, let's start with the proof of the equivalence. Assume first that $$(A_{1}\bowtie H_{1},A_{2}\sharp H_{2})$$ is a Hopf brace in {\sf C}, i.e., the equality 
\begin{align}\label{condHBra_proof}
&\mu_{A\sharp H}^{2}\circ(A\otimes H\otimes \mu_{A\bowtie H}^{1})\\\nonumber=&\mu_{A\bowtie H}^{1}\circ(\mu_{A\sharp H}^{2}\otimes \Gamma_{A_{1}\bowtie H_{1}})\circ(A\otimes H\otimes c_{A\otimes H,A\otimes H}\otimes A\otimes H)\circ (\delta_{A\otimes H}\otimes A\otimes H\otimes A\otimes H)
\end{align}
holds. 

Then, composing on the right of the equality \eqref{condHBra_proof} with $\eta_{A}\otimes H\otimes\eta_{A}\otimes H\otimes A\otimes H$, on the one hand we obtain that
\begin{align*}
\mu_{A\sharp H}^{2}\circ(A\otimes H\otimes \mu_{A\bowtie H}^{1})\circ(\eta_{A}\otimes H\otimes\eta_{A}\otimes H\otimes A\otimes H)=(A\otimes \mu_{H}^{2})\circ(\Psi^{2}\otimes\mu_{H}^{1})\circ(H\otimes\Psi^{1}\otimes H)
\end{align*}
by unit properties and, on the other hand,
\begin{itemize}
\itemindent=-32pt 
\item[ ]$\hspace{0.38cm}\mu_{A\bowtie H}^{1}\circ(\mu_{A\bowtie H}^{2}\otimes \Gamma_{A_{1}\bowtie H_{1}})\circ(A\otimes H\otimes c_{A\otimes H,A\otimes H}\otimes A\otimes H)\circ ((\delta_{A\otimes H}\circ(\eta_{A}\otimes H))\otimes \eta_{A}\otimes H\otimes A\otimes H)$
\item[] $=(A\otimes\mu_{H}^{1})\circ(\Psi^{1}\otimes H)\circ (\mu_{H}^{2}\otimes(\Gamma_{A_{1}\bowtie H_{1}}\circ(\eta_{A}\otimes H\otimes A\otimes H)))\circ(H\otimes c_{H,H}\otimes A\otimes H)\circ(\delta_{H}\otimes H\otimes A\otimes H)$ 
\item[] {\footnotesize\textnormal{(by the condition of coalgebra morphism for $\eta_{A}$, naturality of $c$, \eqref{equal3} and unit properties)}}
\item[] $=(A\otimes \mu_{H}^{1})\circ(((A\otimes\mu_{H}^{1})\circ(\Psi^{1}\otimes H)\circ (H\otimes\Psi^{1}))\otimes\mu_{H}^{2})\circ (\mu_{H}^{2}\otimes ((\lambda_{H}^{1}\otimes \Psi^{2})\circ(\delta_{H}\otimes A))\otimes H)\circ$
\item[] $(H\otimes c_{H,H}\otimes A\otimes H)\circ(\delta_{H}\otimes H\otimes A\otimes H)$ {\footnotesize\textnormal{(by \eqref{equal8} and associativity of $\mu_{H}^{1}$)}}
\item[] $=(A\otimes\mu_{H}^{1})\circ(\Psi^{1}\otimes\mu_{H}^{2})\circ(H\otimes\Psi^{2}\otimes H)\circ(((\Gamma'_{H_{1}}\otimes H)\circ(H\otimes c_{H,H})\circ(\delta_{H}\otimes H))\otimes A\otimes H)$ {\footnotesize\textnormal{(by}}
\item[] {\footnotesize\textnormal{naturality of $c$, coassociativity of $\delta_{H}$ and \eqref{equal5})},}
\end{itemize}
what implies that (C1) holds.

Moreover, composing on the right of the equality \eqref{condHBra_proof} with $A\otimes\eta_{H}\otimes\eta_{A}\otimes H\otimes A\otimes \eta_{H}$, from the left hand hand side the following is obtained:
\begin{align*}
&\mu_{A\sharp H}^{2}\circ(A\otimes H\otimes \mu_{A\bowtie H}^{1})\circ(A\otimes\eta_{H}\otimes\eta_{A}\otimes H\otimes A\otimes \eta_{H})=(\mu_{A}^{2}\otimes H)\circ(A\otimes \Psi^{1})
\end{align*}
by using the unit properties and \eqref{equal2}, while, from the right hand side,
\begin{itemize}
\itemindent=-32pt 
\item[ ]$\hspace{0.38cm}\mu_{A\bowtie H}^{1}\circ(\mu_{A\sharp H}^{2}\otimes \Gamma_{A_{1}\bowtie H_{1}})\circ(A\otimes H\otimes c_{A\otimes H,A\otimes H}\otimes A\otimes H)\circ ((\delta_{A\otimes H}\circ(A\otimes\eta_{H}))\otimes \eta_{A}\otimes H\otimes A\otimes \eta_{H})$
\item[] $=(\mu_{A}^{1}\otimes\mu_{H}^{1})\circ(A\otimes\Psi^{1}\otimes H)\circ(A\otimes H\otimes(\Gamma_{A_{1}\bowtie H_{1}}\circ(A\otimes\eta_{H}\otimes A\otimes\eta_{H})))\circ(A\otimes c_{A,H}\otimes A)\circ(\delta_{A}\otimes H\otimes A)$ 
\item[] {\footnotesize\textnormal{(by the condition of coalgebra morphism for $\eta_{H}$, naturality of $c$, \eqref{equal2} and unit properties)}}
\item [] $=(\mu_{A}^{1}\otimes\mu_{H}^{1})\circ(A\otimes(\Psi^{1}\circ(H\otimes\Gamma_{A_{1}})\circ (c_{A,H}\otimes A))\otimes \eta_{H})\circ(\delta_{A}\otimes H\otimes A)$ {\footnotesize\textnormal{(by \eqref{equal10})}}
\item[] $=(\mu_{A}^{1}\otimes H)\circ(A\otimes(\Psi^{1}\circ(H\otimes\Gamma_{A_{1}})\circ(c_{A,H}\otimes A)))\circ(\delta_{A}\otimes H\otimes A)$ {\footnotesize\textnormal{(by unit property)}.}
\end{itemize}

Then, the equality 
\begin{gather}\label{condC2star}
(\mu_{A}^{2}\otimes H)\circ(A\otimes \Psi^{1})=(\mu_{A}^{1}\otimes H)\circ(A\otimes(\Psi^{1}\circ(H\otimes\Gamma_{A_{1}})\circ(c_{A,H}\otimes A)))\circ(\delta_{A}\otimes H\otimes A)
\end{gather}
holds. As a consequence of the previous equality, (C2) is obtained as follows:
\begin{itemize}
\itemindent=-32pt 
\item[ ]$\hspace{0.38cm}(\Gamma_{A_{1}}\otimes H)\circ(A\otimes\Psi^{1})$
\item[] $=(\mu_{A}^{1}\otimes H)\circ(\lambda_{A}^{1}\otimes ((\mu_{A}^{1}\otimes H)\circ(A\otimes(\Psi^{1}\circ(H\otimes\Gamma_{A_{1}})\circ(c_{A,H}\otimes A)))\circ(\delta_{A}\otimes H\otimes A)))\circ(\delta_{A}\otimes H\otimes A)$ 
\item[] {\footnotesize\textnormal{(by \eqref{condC2star})}}
\item[] $=(\mu_{A}^{1}\otimes H)\circ((\lambda_{A}^{1}\ast id_{A})\otimes (\Psi^{1}\circ(H\otimes\Gamma_{A_{1}})\circ(c_{A,H}\otimes A)))\circ(\delta_{A}\otimes H\otimes A)$ {\footnotesize\textnormal{(by associativity of $\mu_{A}^{1}$ and}}
\item[] {\footnotesize\textnormal{coassociativity of $\delta_{A}$)}}
\item[] $=\Psi^{1}\circ(H\otimes\Gamma_{A_{1}})\circ(c_{A,H}\otimes A)$ {\footnotesize\textnormal{(by \eqref{antipode}, unit property and counit property)}}
\end{itemize}

To conclude this first part of the proof, it only remains to prove that (C3) holds. This condition follows by composing on the right of \eqref{condHBra_proof} with $\eta_{A}\otimes H\otimes A\otimes\eta_{H}\otimes A\otimes\eta_{H}$. On the one hand,
\begin{align*}
\mu_{A\sharp H}^{2}\circ(A\otimes H\otimes \mu_{A\bowtie H}^{1})\circ (\eta_{A}\otimes H\otimes A\otimes\eta_{H}\otimes A\otimes\eta_{H})=\Psi^{2}\circ(H\otimes\mu_{A}^{1})
\end{align*}
by unit properties and \eqref{equal1}. On the other hand,
\begin{itemize}
\itemindent=-32pt 
\item[ ]$\hspace{0.38cm}\mu_{A\bowtie H}^{1}\circ(\mu_{A\sharp H}^{2}\otimes \Gamma_{A_{1}\bowtie H_{1}})\circ(A\otimes H\otimes c_{A\otimes H,A\otimes H}\otimes A\otimes H)\circ ((\delta_{A\otimes H}\circ(\eta_{A}\otimes H))\otimes A\otimes \eta_{H}\otimes A\otimes \eta_{H})$
\item [] $=(\mu_{A}^{1}\otimes\mu_{H}^{1})\circ(A\otimes\Psi^{1}\otimes H)\circ(\Psi^{2}\otimes(\Gamma_{A_{1}\bowtie H_{1}}\circ (\eta_{A}\otimes H\otimes A\otimes\eta_{H})))\circ(H\otimes c_{H,A}\otimes A)\circ(\delta_{H}\otimes A\otimes A)$ 
\item[] {\footnotesize\textnormal{(by the condition of coalgebra morphism for $\eta_{A}$, naturality of $c$ and unit properties)}}
\item [] $=(\mu_{A}^{1}\otimes\mu_{H}^{1})\circ(A\otimes\Psi^{1}\otimes H)\circ(\Psi^{2}\otimes((A\otimes\mu_{H}^{1})\circ (\Psi^{1}\otimes H)\circ (\lambda_{H}^{1}\otimes\Psi^{2})\circ (\delta_{H}\otimes A)))\circ(H\otimes c_{H,A}\otimes $
\item[] $A)\circ(\delta_{H}\otimes A\otimes A)$ {\footnotesize\textnormal{(by \eqref{equal9})}}
\item [] $=(\mu_{A}^{1}\otimes\mu_{H}^{1})\circ(A\otimes((A\otimes\mu_{H}^{1})\circ(\Psi^{1}\otimes H)\circ(H\otimes\Psi^{1}))\otimes H)\circ(A\otimes H\otimes ((\lambda_{H}^{1}\otimes \Psi^{2})\circ(\delta_{H}\otimes A)))\circ$
\item[] $(((\Psi^{2}\otimes H)\circ(H\otimes c_{H,A})\circ(\delta_{H}\otimes A))\otimes A)$ {\footnotesize\textnormal{(by naturality of $c$ and associativity of $\mu_{H}^{1}$)}}
\item[] $=(\mu_{A}^{1}\otimes\mu_{H}^{1})\circ(A\otimes((\Psi^{1}\otimes H)\circ((id_{H}\ast\lambda_{H}^{1})\otimes\Psi^{2})\circ(\delta_{H}\otimes A)))\circ(\Psi^{2}\otimes A)$ {\footnotesize\textnormal{(by \eqref{equal5}, \eqref{equal6} and coassociativity}}
\item[] {\footnotesize\textnormal{of $\delta_{H}$)}}
\item[] $=(\mu_{A}^{1}\otimes H)\circ (A\otimes\Psi^{2})\circ(\Psi^{2}\otimes A)$ {\footnotesize\textnormal{(by \eqref{equal1} and (co)unit properties)},}
\end{itemize}
and, as a consequence, (C3) is satisfied.

Let's prove the remaining implication. If conditions (C1), (C2) and (C3) are assumed, then we obtain that
\begin{itemize}
\itemindent=-32pt 
\item[ ]$\hspace{0.38cm}\mu_{A\bowtie H}^{1}\circ(\mu_{A\sharp H}^{2}\otimes \Gamma_{A_{1}\bowtie H_{1}})\circ(A\otimes H\otimes c_{A\otimes H,A\otimes H}\otimes A\otimes H)\circ (\delta_{A\otimes H}\otimes A\otimes H\otimes A\otimes H)$
\item[] $=(\mu_{A}^{1}\otimes\mu_{H}^{1})\circ(A\otimes\Psi^{1}\otimes H)\circ(((\mu_{A}^{2}\otimes\mu_{H}^{2})\circ(A\otimes\Psi^{2}\otimes H))\otimes((A\otimes\mu_{H}^{1})\circ((\Psi^{1}\circ (\lambda_{H}^{1}\otimes\Gamma_{A_{1}}))\otimes\mu_{H}^{2})\circ$
\item[] $(((c_{A,H}\otimes\Psi^{2})\circ (A\otimes\delta_{H}\otimes A))\otimes H)))\circ (A\otimes H\otimes c_{A\otimes H,A\otimes H}\otimes A\otimes H)\circ(\delta_{A\otimes H}\otimes A\otimes H\otimes A\otimes H)$ 
\item[] {\footnotesize\textnormal{(by \eqref{equal7})}}
\item[] $=(\mu_{A}^{1}\otimes\mu_{H}^{1})\circ (\mu_{A}^{2}\otimes((A\otimes\mu_{H}^{1})\circ (\Psi^{1}\otimes H)\circ(H\otimes\Psi^{1})\circ (\mu_{H}^{2}\otimes((\lambda_{H}^{1}\otimes\Gamma_{A_{1}})\circ (c_{A,H}\otimes A))))\otimes\mu_{H}^{2})\circ $
\item[] $(A\otimes((A\otimes H\otimes c_{A,H})\circ (A\otimes c_{A,H}\otimes H)\circ (c_{A,A}\otimes H\otimes H))\otimes((H\otimes A\otimes\mu_{H}^{2})\circ (H\otimes\Psi^{2}\otimes H)\circ(\delta_{H}\otimes$
\item[] $A\otimes H)))\circ(\delta_{A}\otimes ((\Psi^{2}\otimes c_{H,H})\circ (H\otimes c_{H,A}\otimes H)\circ (\delta_{H}\otimes A\otimes H))\otimes A\otimes H)$ {\footnotesize\textnormal{(by naturality of $c$ and}}
\item[] {\footnotesize\textnormal{associativity of $\mu_{H}^{1}$)}}
\item[] $=(\mu_{A}^{1}\otimes\mu_{H}^{1})\circ(\mu_{A}^{2}\otimes(\Psi^{1}\circ (H\otimes\Gamma_{A_{1}})\circ(c_{A,H}\otimes A))\otimes\mu_{H}^{2})\circ (A\otimes c_{A,A}\otimes((\Gamma'_{H_{1}}\otimes\Psi^{2})\circ(H\otimes c_{H,H}\otimes $
\item[] $A)\circ(\delta_{H}\otimes H\otimes A))\otimes H)\circ(\delta_{A}\otimes\Psi^{2}\otimes H\otimes A\otimes H)$ {\footnotesize\textnormal{(by \eqref{equal5}, \eqref{equal6}, coassociativity of $\delta_{H}$ and naturality of $c$)}}
\item[] $=((\mu_{A}^{1}\circ(\mu_{A}^{2}\otimes\Gamma_{A_{1}})\circ(A\otimes c_{A,A}\otimes A)\circ(\delta_{A}\otimes A\otimes A))\otimes H)\circ(A\otimes A\otimes((A\otimes\mu_{H}^{1})\circ(\Psi^{1}\otimes\mu_{H}^{2})\circ$
\item[] $(H\otimes\Psi^{2}\otimes H)\circ (((\Gamma'_{H_{1}}\otimes H)\circ(H\otimes c_{H,H})\circ(\delta_{H}\otimes H))\otimes A\otimes H)))\circ(A\otimes\Psi^{2}\otimes H\otimes A\otimes H)$ {\footnotesize\textnormal{(by (C2))}}
\item[] $=(\mu_{A}^{2}\otimes\mu_{H}^{2})\circ (A\otimes((\mu_{A}^{1}\otimes H)\circ (A\otimes\Psi^{2})\circ(\Psi^{2}\otimes A))\otimes\mu_{H}^{1})\circ(A\otimes H\otimes A\otimes\Psi^{1}\otimes H)$ {\footnotesize\textnormal{(by (C1) and \eqref{compatHbrace}}}
\item[] {\footnotesize\textnormal{for $\mathbb{A}$)}}
\item[] $=\mu_{A\sharp H}^{2}\circ(A\otimes H\otimes\mu_{A\bowtie H}^{1})$ {\footnotesize\textnormal{(by (C3))}}
\end{itemize}
what concludes the proof.
\end{proof}

As a corollary of the previous theorem, A. Agore's Theorem \ref{agtheorem2} is obtained.
\begin{corollary}\label{cor_Agore1}
Let $A$ be a Hopf algebra and $\mathbb{H}=(H_{1},H_{2})$ a cocommutative Hopf brace in {\sf C}. If the conditions
\begin{itemize}
\item[(i)] $(A,H_{1},\varphi_{A}^{1},\phi_{H})$ is a matched pair of Hopf algebras,
\item[(ii)] $(A,\varphi_{A}^{2})$ is a left $H_{2}$-module algebra-coalgebra,
\item[(iii)] $\varphi_{A}^{2}\circ(H\otimes\varphi_{A}^{1})=\varphi_{A}^{1}\circ (\Gamma'_{H_{1}}\otimes \varphi_{A}^{2})\circ (H\otimes c_{H,H}\otimes A)\circ(\delta_{H}\otimes H\otimes A)$,
\item[(iv)] $\mu_{H}^{2}\circ(H\otimes\phi_{H})=\mu_{H}^{1}\circ(\phi_{H}\otimes H)\circ (\Gamma'_{H_{1}}\otimes\Psi^{2})\circ (H\otimes c_{H,H}\otimes A)\circ (\delta_{H}\otimes H\otimes A)$
\end{itemize}
hold, then
$(A\bowtie H_{1},A\sharp H_{2})$
is a Hopf brace in {\sf C}.
\end{corollary}
\begin{proof}
It is enough to apply Theorem \ref{mainth} considering the trivial Hopf brace obtained from the Hopf algebra $A$, $\mathbb{A}=\mathbb{A}_{triv}=(A,A)$. Note that, due to the cocommutativity of $\mathbb{H}$, \eqref{adcond_mp_phiHtriv} always holds in this situation. So, let's compute conditions (C1), (C2) and (C3). At first, by associativity of $\mu_{A}$, \eqref{antipode} and unit property, it results straightforward to prove that $$\Gamma_{A}=\varepsilon_{A}\otimes A$$ and, as a consequence, (C2) always holds. Moreover,
\begin{itemize}
\itemindent=-32pt 
\item[ ]$\hspace{0.38cm}(\mu_{A}\otimes H)\circ(A\otimes\Psi^{2})\circ(\Psi^{2}\otimes A)$
\item[] $=((\mu_{A}\circ(\varphi_{A}^{2}\otimes\varphi_{A}^{2})\circ (H\otimes c_{H,A}\otimes A)\circ(\delta_{H}\otimes A\otimes A))\otimes H)\circ(H\otimes((A\otimes c_{H,A})\circ (c_{H,A}\otimes A)))\circ(\delta_{H}\otimes A\otimes A)$ 
\item[] {\footnotesize\textnormal{(by naturality of $c$ and coassociativity of $\delta_{H}$)}}
\item[] $=\Psi^{2}\circ (H\otimes\mu_{A})$ {\footnotesize\textnormal{(by the condition of morphism of left $H_{2}$-modules for $\mu_{A}$ and naturality of $c$)},}
\end{itemize}
what implies that (C3) always holds too in this situation. To conclude, let's see that hypothesis (iii) and (iv) imply condition (C1). First of all, note that the following equalities hold:
\begin{gather}
\label{condproof1}
(\delta_{A}\otimes H)\circ\Psi^{2}=(A\otimes\Psi^{2})\circ(\Psi^{2}\otimes A)\circ(H\otimes\delta_{A}),
\end{gather}
what follows by the condition of coalgebra morphism for $\varphi_{A}^{2}$, naturality of $c$ and coassociativity of $\delta_{H}$, and
\begin{gather}\label{condproof2}
\Gamma'_{H_{1}}\circ(H\otimes\phi_{H})=\phi_{H}\circ(\Gamma'_{H_{1}}\otimes\varphi_{A}^{2})\circ(H\otimes c_{H,H}\otimes A)\circ(\delta_{H}\otimes H\otimes A),
\end{gather}
which is a consequence of hypothesis (iv) as we can see in what follows:
\begin{itemize}
\itemindent=-32pt 
\item[ ]$\hspace{0.38cm}\Gamma'_{H_{1}}\circ(H\otimes\phi_{H})$
\item [] $=\mu_{H}^{1}\circ ((\mu_{H}^{2}\circ(H\otimes\phi_{H}))\otimes\lambda_{H}^{1})\circ(H\otimes((H\otimes c_{H,A})\circ (c_{H,H}\otimes A)))\circ(\delta_{H}\otimes H\otimes A)$ {\footnotesize\textnormal{(by naturality of $c$)}}
\item[] $=\mu_{H}^{1}\circ((\mu_{H}^{1}\circ(\phi_{H}\otimes H)\circ (\Gamma'_{H_{1}}\otimes\Psi^{2})\circ (H\otimes c_{H,H}\otimes A)\circ (\delta_{H}\otimes H\otimes A))\otimes\lambda_{H}^{1})\circ(H\otimes((H\otimes c_{H,A})\circ $
\item[] $(c_{H,H}\otimes A)))\circ(\delta_{H}\otimes H\otimes A)$ {\footnotesize\textnormal{(by (iv))}}
\item[] $=\mu_{H}^{1}\circ (\phi_{H}\otimes H)\circ(\Gamma'_{H_{1}}\otimes ((A\otimes\mu_{H}^{1})\circ (\Psi^{2}\otimes\lambda_{H}^{1})\circ (H\otimes c_{H,A})\circ (\delta_{H}\otimes A)))\circ(H\otimes c_{H,H}\otimes A)\circ(\delta_{H}\otimes H\otimes A)$ 
\item[]{\footnotesize\textnormal{(by coassociativity of $\delta_{H}$, naturality of $c$ and asocciativity of $\mu_{H}^{1}$)}}
\item[] $=\mu_{H}^{1}\circ(\phi_{H}\otimes(id_{H}\ast\lambda_{H}^{1}))\circ (\Gamma'_{H_{1}}\otimes\Psi^{2})\circ(H\otimes c_{H,H}\otimes A)\circ(\delta_{H}\otimes H\otimes A)$ {\footnotesize\textnormal{(by \eqref{equal6})}}
\item[] $=\phi_{H}\circ(\Gamma'_{H_{1}}\otimes \varphi_{A}^{2})\circ (H\otimes c_{H,H}\otimes A)\circ(\delta_{H}\otimes H\otimes A)$ {\footnotesize\textnormal{(by \eqref{antipode}, unit property and \eqref{equaleps})}.}
\end{itemize}

Therefore, we obtain:
\begin{itemize}
\itemindent=-32pt 
\item[ ]$\hspace{0.38cm}(A\otimes\mu_{H}^{1})\circ(\Psi^{1}\otimes\mu_{H}^{2})\circ(H\otimes\Psi^{2}\otimes H)\circ(((\Gamma'_{H_{1}}\otimes H)\circ(H\otimes c_{H,H})\circ(\delta_{H}\otimes H))\otimes A\otimes H)$
\item[] $=(\varphi_{A}^{1}\otimes (\mu_{H}^{1}\circ (\phi_{H}\otimes\mu_{H}^{2})))\circ(H\otimes c_{H,A}\otimes A\otimes H\otimes H)\circ(((\Gamma'_{H_{1}}\otimes\Gamma'_{H_{1}})\circ(H\otimes c_{H,H}\otimes H)\circ(\delta_{H}\otimes$
\item[] $\delta_{H}))\otimes((A\otimes\Psi^{2})\circ (\Psi^{2}\otimes A)\circ (H\otimes\delta_{A}))\otimes H)\circ(H\otimes c_{H,H}\otimes A\otimes H)\circ(\delta_{H}\otimes H\otimes A\otimes H)$ {\footnotesize\textnormal{(by the}}
\item[] {\footnotesize\textnormal{condition of coalgebra morphism for $\Gamma'_{H_{1}}$ and \eqref{condproof1})}}
\item []$=((\varphi_{A}^{1}\circ(\Gamma'_{H_{1}}\otimes\varphi_{A}^{2})\circ (H\otimes c_{H,H}\otimes A)\circ(\delta_{H}\otimes H\otimes A))\otimes(\mu_{H}^{1}\circ (\phi_{H}\otimes\mu_{H}^{2})\circ (\Gamma'_{H_{1}}\otimes\Psi^{2}\otimes H)\circ (H\otimes $
\item [] $c_{H,H}\otimes A\otimes H)))\circ (H\otimes (c_{H\otimes H,H\otimes A}\circ(\delta_{H}\otimes H\otimes A))\otimes H\otimes A\otimes H)\circ (\delta_{H}\otimes\delta_{H\otimes A}\otimes H)$ {\footnotesize\textnormal{(by naturality}}
\item[] {\footnotesize\textnormal{of $c$ and cocommutativity of $\delta_{H}$)}}
\item[] $=(A\otimes (\mu_{H}^{1}\circ (\phi_{H}\otimes\mu_{H}^{2})\circ (\Gamma'_{H_{1}}\otimes\Psi^{2}\otimes H)\circ (H\otimes c_{H,H}\otimes A\otimes H)\circ (\delta_{H}\otimes H\otimes A\otimes H)))\circ(\Psi^{2}\otimes H\otimes $
\item[] $A\otimes H)\circ (H\otimes ((\varphi_{A}^{1}\otimes H\otimes A)\circ (H\otimes c_{H,A}\otimes A)\circ(\delta_{H}\otimes\delta_{A}))\otimes H)$ {\footnotesize\textnormal{(by (iii) and naturality of $c$)}}
\item[] $=(A\otimes(\mu_{H}^{1}\circ ((\phi_{H}\circ (\Gamma'_{H_{1}}\otimes\varphi_{A}^{2})\circ (H\otimes c_{H,H}\otimes A)\circ(\delta_{H}\otimes H\otimes A))\otimes\mu_{H}^{2})\circ (H\otimes((H\otimes c_{H,A})\circ (c_{H,H}\otimes A))\otimes $
\item[] $H)\circ (\delta_{H}\otimes H\otimes A\otimes H)))\circ(\Psi^{2}\otimes H\otimes A\otimes H)\circ (H\otimes ((\varphi_{A}^{1}\otimes H\otimes A)\circ (H\otimes c_{H,A}\otimes A)\circ(\delta_{H}\otimes\delta_{A}))\otimes H)$ 
\item[] {\footnotesize\textnormal{(by naturality of $c$ and coassociativity of $\delta_{H}$)}}
\item[] $=(A\otimes (\mu_{H}^{1}\circ (\Gamma'_{H_{1}}\otimes\mu_{H}^{2})\circ (H\otimes c_{H,H}\otimes H)\circ (\delta_{H}\otimes H\otimes H)))\circ(\Psi^{2}\otimes H\otimes H)\circ(H\otimes\Psi^{1}\otimes H)$ {\footnotesize\textnormal{(by}}
\item[] {\footnotesize\textnormal{\eqref{condproof2} and naturality of $c$)}}
\item[] $=(A\otimes \mu_{H}^{2})\circ(\Psi^{2}\otimes\mu_{H}^{1})\circ(H\otimes\Psi^{1}\otimes H)$ {\footnotesize\textnormal{(by \eqref{compatHbrace2} for $\mathbb{H}$)},}
\end{itemize}
and then (C1) is satisfied. In conclusion, if (C1), (C2) and (C3) hold, then $(A\bowtie H_{1},A\sharp H_{2})$ is a Hopf brace in {\sf C} by Theorem \ref{mainth}.
\end{proof}
\begin{remark}
Assume the conditions of the previous corollary. In the proof we have seen that if (iii) and (iv) hold, then (C1) is verified. However, this result is in fact an equivalence. Let's see that condition (C1) implies that (iii) and (iv) hold too. Indeed, composing on the right of equality (C1) with $H\otimes H\otimes A\otimes\eta_{H}$ and on the left with $A\otimes\varepsilon_{H}$, (iii) is obtained by using the condition of coalgebra morphism for $\mu_{H}^{1}$, $\mu_{H}^{2}$ and $\phi_{H}$ and the counit properties. On the other hand, (iv) is obtained by composing on the right of equality (C1) with $H\otimes H\otimes A\otimes\eta_{H}$, on the left with $\varepsilon_{A}\otimes H$ and by using the condition of coalgebra morphism for $\varphi_{A}^{1}$ and $\varphi_{A}^{2}$ and the (co)unit properties.
\end{remark}
\end{remark}

\begin{corollary}\label{smash_cor_mainth}
Let $A$ and $H$ be Hopf algebras in {\sf C} such that 
\begin{itemize}
\item[(i)] $(A,H,\varphi_{A}^{1},\phi_{H})$ is a matched pair of Hopf algebras in {\sf C},
\item[(ii)] $(A,\varphi_{A}^{2})$ is a left $H$-module algebra-coalgebra satisfying \eqref{adcond_mp_phiHtriv}.
\end{itemize}

Then, the pair $(A\bowtie H,A\sharp H)$ is a Hopf brace in {\sf C} if and only if equality
\begin{equation}\label{C1_smash}
(A\otimes \mu_{H})\circ(\Psi^{1}\otimes H)\circ(\varphi_{H}^{ad}\otimes\Psi^{2})\circ(H\otimes c_{H,H}\otimes A)\circ(\delta_{H}\otimes H\otimes A)=(A\otimes\mu_{H})\circ(\Psi^{2}\otimes H)\circ(H\otimes \Psi^{1}),
\end{equation}
where $\Psi^{1}$ and $\Psi^{2}$ are defined as in the statement of Theorem \ref{mainth}.
\end{corollary}
\begin{proof}
It is enough to apply Theorem \ref{mainth} by considering that $\mathbb{A}=\mathbb{A}_{triv}=(A,A)$ and $\mathbb{H}=\mathbb{H}_{triv}=(H,H)$. As was proved in Corollary \ref{cor_Agore1}, when $\mathbb{A}=\mathbb{A}_{triv}$, conditions (C2) and (C3) of Theorem \ref{mainth} automatically hold. So, it is enough to compute (C1) assuming that $\mathbb{H}=\mathbb{H}_{triv}$. In this situation, note first that $\Gamma'_{H}=\varphi_{H}^{ad}$, then (C1) is equivalent to 
\begin{align}\label{t1}
&(A\otimes\mu_{H})\circ(((A\otimes\mu_{H})\circ(\Psi^{1}\otimes H)\circ(\varphi_{H}^{ad}\otimes\Psi^{2})\circ(H\otimes c_{H,H}\otimes A)\circ(\delta_{H}\otimes H\otimes A))\otimes H)\\\nonumber=&(A\otimes\mu_{H})\circ(((A\otimes\mu_{H})\circ(\Psi^{2}\otimes H)\circ(H\otimes\Psi^{1}))\otimes H)
\end{align}
by associativity of $\mu_{H}$. Then, composing on the right with $H\otimes H\otimes A\otimes\eta_{H}$, it is obtained that \eqref{t1} is equivalent to \eqref{C1_smash}.
\end{proof}

\begin{corollary}\label{tensor_cor_mainth}
Let $A$ and $H$ be Hopf algebras in {\sf C} such that $(A,H,\varphi_{A},\phi_{H})$ is a matched pair. Then, the pair $(A\bowtie H,A\otimes H)$ is a Hopf brace if and only if 
\begin{equation}\label{C1_tensor}
(A\otimes \mu_{H})\circ(\Psi\otimes H)\circ(\varphi_{H}^{ad}\otimes c_{H,A})\circ(H\otimes c_{H,H}\otimes A)\circ(\delta_{H}\otimes H\otimes A)=(A\otimes\mu_{H})\circ(c_{H,A}\otimes H)\circ(H\otimes\Psi),
\end{equation}
where $\Psi$ is defined as in \eqref{interw}.
\end{corollary}
\begin{proof}
This result is a direct consequence of Corollary \ref{smash_cor_mainth} by considering that $\varphi_{A}^{2}=\varepsilon_{H}\otimes A$. Under this assumption, \eqref{adcond_mp_phiHtriv} automatically holds and $\Psi^{2}=c_{H,A}$, what implies that \eqref{C1_smash} becomes \eqref{C1_tensor} and also $A\sharp H=A\otimes H$ (see Remark \ref{particular_smash}).
\end{proof}

\begin{corollary}\label{tensor_cor_mainth_commut}
Let $A$ and $H$ be Hopf algebras in {\sf C} such that $(A,H,\varphi_{A},\phi_{H})$ is a matched pair. If $H$ is commutative, then $(A\bowtie H,A\otimes H)$ is a Hopf brace in {\sf C}.
\end{corollary}
\begin{proof}
This result is a consequence of Corollary \ref{tensor_cor_mainth} because \eqref{C1_tensor} always holds when $H$ is commutative. Indeed, under commutativity of $H$, $\varphi_{H}^{ad}=\varepsilon_{H}\otimes H$ and then we have the following:
\begin{itemize}
\itemindent=-32pt 
\item[ ]$\hspace{0.38cm}(A\otimes \mu_{H})\circ(\Psi\otimes H)\circ(\varphi_{H}^{ad}\otimes c_{H,A})\circ(H\otimes c_{H,H}\otimes A)\circ(\delta_{H}\otimes H\otimes A)$
\item[] $=(A\otimes \mu_{H})\circ(\Psi\otimes H)\circ(H\otimes c_{H,A})\circ(\varepsilon_{H}\otimes c_{H,H}\otimes A)\circ(\delta_{H}\otimes H\otimes A)$ {\footnotesize\textnormal{(by $\varphi_{H}^{ad}=\varepsilon_{H}\otimes H$)}}
\item[] $=(A\otimes\mu_{H})\circ(\Psi\otimes H)\circ(H\otimes c_{H,A})\circ(c_{H,H}\otimes A)$ {\footnotesize\textnormal{(by counit property)}}
\item[] $=(A\otimes (\mu_{H}\circ c_{H,H}))\circ(c_{H,A}\otimes H)\circ(H\otimes \Psi)$ {\footnotesize\textnormal{(by naturality of $c$)}}
\item[] $=(A\otimes\mu_{H})\circ(c_{H,A}\otimes H)\circ(H\otimes \Psi)$ {\footnotesize\textnormal{(by commutativity of $H$).}}\qedhere
\end{itemize}
\end{proof}

\begin{corollary}\label{D(H)_commutative_HBr}
Let $H$ be a finite and commutative Hopf algebra in {\sf C}. Then, $(D(H),T(H))$ is a Hopf brace in {\sf C}.
\end{corollary}
\begin{proof}
This result follows directly from Corollary \ref{tensor_cor_mainth_commut} by considering the particular matched pair $(\wh,H,\varphi_{\wh},$ $\phi_{H})$, whose actions involved are defined in \eqref{varphi_wh} and \eqref{phi_H}, respectively.
\end{proof}

\begin{example} Assume that ${\sf C}={}_{\mathbb{F}}{\sf Vect}.$ By Example \ref{D(F[G])_cocom} it is already known that, if $(G,\cdot)$ is a finite group, $\mathbb{F}[G]$ is a finite and cocommutative Hopf algebra in ${}_{\mathbb{F}}{\sf Vect}$ whose dual Hopf algebra $\mathbb{F}[G]^{\ast}$ is commutative. Then, by Corollary \ref{D(H)_commutative_HBr}, $(D(\mathbb{F}[G]^{\ast}),T(\mathbb{F}[G]^{\ast}))$ is a Hopf brace in ${}_{\mathbb{F}}{\sf Vect}$. 

Moreover, if $(G,\cdot)$ is commutative, then $\mathbb{F}[G]$ is also a commutative Hopf algebra. Thus, again by Corollary \ref{D(H)_commutative_HBr}, it is obtained that $(D(\mathbb{F}[G]),T(\mathbb{F}[G]))$ is a Hopf brace in ${}_{\mathbb{F}}{\sf Vect}$.
\end{example}

However, to obtain conditions under which $(D(H),T(H))$ is a Hopf brace when $H$ is not commutative we have to resort to Corollary \ref{tensor_cor_mainth}.
\begin{corollary}\label{D(H)_noncomm_HBr}
Let $H$ be a finite Hopf algebra in {\sf C}. The pair $(D(H),T(H))$ is a Hopf brace in {\sf C} if and only if 
\begin{align}\label{D(H)_HBr_nocommutativity}
&(\mu_{H}\otimes H)\circ(H\otimes c_{H,H})\circ(J\otimes H)\circ(H\otimes((\varphi_{H}^{ad}\otimes H)\circ(H\otimes c_{H,H})\circ(\delta_{H}\otimes H)))\\\nonumber=&(\mu_{H}\otimes H)\circ(H\otimes J)\circ(c_{H,H}\otimes H),
\end{align}
where $J$ is defined by 
\begin{align}\label{J_def}
J\coloneqq (H\otimes(\mu_{H}\circ(\lambda_{H}^{-1}\otimes H)))\circ(\delta_{H}\otimes H)\circ c_{H,H}\circ(\mu_{H}\otimes H)\circ(H\otimes\delta_{H}).
\end{align}
\end{corollary}
\begin{proof}
It is enough to apply Corollary \ref{tensor_cor_mainth} by considering the matched pair $(\wh,H,\varphi_{\wh},$ $\phi_{H})$, whose actions involved are defined in \eqref{varphi_wh} and \eqref{phi_H}, respectively, from which Drinfeld's Double is constructed. Therefore, we have to particularize \eqref{C1_tensor} to this situation.

First of all, let's compute $\Psi$ for $D(H)=\wh\bowtie H.$ Using \eqref{psi_eta} and \eqref{for_prod_D(H)}, it is obtained that the required $\Psi$ admits the following expression:
\begin{equation}\label{psi_D(H)}
\Psi=(b_{H}(K)\otimes c_{H, H^{\ast}})\circ(\delta_{H\otimes\wh}\otimes(b_{H}(K)\circ(\lambda_{H}^{-1}\otimes H^{\ast})))\circ\delta_{H\otimes\wh},
\end{equation}
which can be expanded obtaining that
\begin{align}\label{psi_expanded}
\Psi=(H^{\ast}\otimes H\otimes b_{H}(K))\circ(H^{\ast}\otimes J\otimes H^{\ast})\circ(a_{H}(K)\otimes H\otimes H^{\ast}).
\end{align}

Indeed: Note that 
\begin{align}\label{x1}
&(H^{\ast}\otimes b_{H}(K))\circ (c_{H,H^{\ast}}\otimes H^{\ast})\circ(H\otimes\delta_{\wh})\\\nonumber=&(((H^{\ast}\otimes b_{H}(K))\circ (c_{H,H^{\ast}}\otimes H^{\ast}))\otimes(b_{H}(K)\circ (\mu_{H}\otimes H^{\ast})))\circ(H\otimes((H^{\ast}\otimes a_{H}(K)\otimes H)\\\nonumber&\hspace{-.38 cm}\circ a_{H}(K))\otimes H^{\ast})\;\footnotesize\textnormal{(by definition of $\delta_{\wh}$)}\\\nonumber=&(H^{\ast}\otimes(b_{H}(K)\circ (\mu_{H}\otimes H^{\ast})))\circ(((c_{H,H^{\ast}}\otimes H)\circ(H\otimes a_{H}(K)))\otimes H^{\ast})\;\footnotesize\textnormal{(by the properties of}\\\nonumber& \hspace{-.30cm}\footnotesize\textnormal{the unit and the counit of the adjunction)}\\\nonumber=&(H^{\ast}\otimes b_{H}(K))\circ(H^{\ast}\otimes(\mu_{H}\circ c_{H,H})\otimes H^{\ast})\circ(a_{H}(K)\otimes H\otimes H^{\ast})\;\footnotesize\textnormal{(by \eqref{fin_unit})}
\end{align}
and, using similar arguments, it also holds that
\begin{align}\label{x2}
(b_{H}(K)\otimes H^{\ast})\circ(H\otimes\delta_{\wh})=(H^{\ast}\otimes b_{H}(K))\circ(H^{\ast}\otimes\mu_{H}\otimes H^{\ast})\circ(a_{H}(K)\otimes H\otimes H^{\ast}).
\end{align}

Then, by \eqref{psi_D(H)}, \eqref{psi_expanded} follows by
\begin{itemize}
\itemindent=-32pt 
\item[ ]$\hspace{0.38cm}\Psi$
\item[] $=(((b_{H}(K)\otimes H^{\ast})\circ(H\otimes\delta_{\wh}))\otimes H)\circ (H\otimes c_{H,H^{\ast}})\circ(H\otimes H\otimes ((H^{\ast}\otimes b_{H}(K))\circ (c_{H,H^{\ast}}\otimes H^{\ast})\circ$
\item[] $(H\otimes\delta_{\wh})))\circ(((\delta_{H}\otimes\lambda_{H}^{-1})\circ\delta_{H})\otimes H^{\ast})$ {\footnotesize\textnormal{(by naturality of $c$)}}
\item[] $=(H^{\ast}\otimes ((b_{H}(K)\otimes H)\circ(\mu_{H}\otimes c_{H,H^{\ast}}))\otimes (b_{H}(K)\circ((\mu_{H}\circ c_{H,H})\otimes H^{\ast})))\circ(a_{H}(K)\otimes ((\delta_{H}\otimes a_{H}(K)\otimes$
\item[] $\lambda_{H}^{-1})\circ\delta_{H})\otimes H^{\ast})$ {\footnotesize\textnormal{(by \eqref{x1} and \eqref{x2})}}
\item[] $=(H^{\ast}\otimes ((H\otimes b_{H}(K))\circ((c_{H,H}\circ (\mu_{H}\otimes H))\otimes H^{\ast}))\otimes (b_{H}(K)\circ ((\mu_{H}\circ c_{H,H})\otimes H^{\ast})))\circ(a_{H}(K)\otimes$
\item[] $ ((\delta_{H}\otimes a_{H}(K)\otimes\lambda_{H}^{-1})\circ\delta_{H})\otimes H^{\ast})$ {\footnotesize\textnormal{(by \eqref{fin_counit})}}
\item[] $=(H^{\ast}\otimes H\otimes b_{H}(K))\circ(H^{\ast}\otimes((H\otimes(\mu_{H}\circ c_{H,H}))\circ((c_{H,H}\circ (\mu_{H}\otimes H)\circ (H\otimes\delta_{H}))\otimes\lambda_{H}^{-1})\circ(H\circ\delta_{H}))\otimes H^{\ast})\circ$
\item[] $(a_{H}(K)\otimes H\otimes H^{\ast})$ {\footnotesize\textnormal{(by the properties of the unit and the counit of the adjunction)}}
\item[] $=(H^{\ast}\otimes H\otimes b_{H}(K))\circ(H^{\ast}\otimes J\otimes H^{\ast})\circ(a_{H}(K)\otimes H\otimes H^{\ast})$ {\footnotesize\textnormal{(by coassociativity of $\delta_{H}$, naturality of $c$ and }}
\item[] {\footnotesize\textnormal{\eqref{J_def})}}.
\end{itemize}

In this case, by \eqref{psi_expanded}, \eqref{C1_tensor} holds if and only if 
\begin{align}\label{b1}
&(H^{\ast}\otimes\mu_{H}\otimes b_{H}(K))\circ(H^{\ast}\otimes H\otimes J\otimes H^{\ast})\circ(((c_{H,H^{\ast}}\otimes H)\circ(H\otimes a_{H}(K)))\otimes H\otimes H^{\ast})\\\nonumber=&(H^{\ast}\otimes\mu_{H})\circ (H^{\ast}\otimes((H\otimes b_{H}(K))\circ (J\otimes H^{\ast}))\otimes H)\circ(a_{H}(K)\otimes((\varphi_{H}^{ad}\otimes c_{H,H^{\ast}})\circ (H\otimes \\\nonumber&\hspace{-.30 cm}c_{H,H}\otimes H^{\ast})\circ (\delta_{H}\otimes H\otimes H^{\ast}))).
\end{align}

Furthermore, by \eqref{fin_unit} and \eqref{fin_counit}, \eqref{b1} can be rewritten as follows:
\begin{align*}
&(H^{\ast}\otimes H\otimes b_{H}(K))\circ(H^{\ast}\otimes((\mu_{H}\otimes H)\circ(H\otimes J)\circ(c_{H,H}\otimes H))\otimes H^{\ast})\circ(a_{H}(K)\otimes H\otimes H\otimes H^{\ast})\\=&(H^{\ast}\otimes H\otimes b_{H}(K))\circ(H^{\ast}\otimes((\mu_{H}\otimes H)\circ(H\otimes c_{H,H})\circ(J\otimes H))\otimes H^{\ast})\circ(H^{\ast}\otimes H\otimes((\varphi_{H}^{ad}\otimes H)\circ\\&\hspace{-.3 cm}(H\otimes c_{H,H})\circ(\delta_{H}\otimes H))\otimes H^{\ast})\circ(a_{H}(K)\otimes H\otimes H\otimes H^{\ast})
\end{align*}
which is equivalent to \eqref{D(H)_HBr_nocommutativity} by the properties of the unit and the counit of the {\sf C}-adjunction $H\otimes - \dashv H^{\ast}\otimes -$.
\end{proof}

\section*{Authors Contribution} Ram\'on Gonz\'alez Rodr\'{\i}guez and Brais Ramos Pérez  participated at all stages in the preparation of the manuscript.

\section*{Funding Declaration}
The  authors were supported by  Ministerio de Ciencia e Innovaci\'on of Spain. Agencia Estatal de Investigaci\'on. Uni\'on Europea - Fondo Europeo de Desarrollo Regional (FEDER). Grant PID2020-115155GB-I00: Homolog\'{\i}a, homotop\'{\i}a e invariantes categ\'oricos en grupos y \'algebras no asociativas.

Moreover, Brais Ramos Pérez was funded by Xunta de Galicia, grant ED431C 2023/31 (European FEDER support included, UE).

Also, Brais Ramos Pérez was financially supported by Xunta de Galicia Scholarship ED481A-2023-023.

\section*{Data Availability} Not applicable to this article.

\section*{Competing Interests Declaration} The authors declares no conflict of interest.

\bibliographystyle{amsalpha}

\end{document}